\documentclass[preprint,11pt,1p,nopreprintline]{elsarticle}

\usepackage{graphicx}
\usepackage[all,cmtip]{xy}
\usepackage{amsthm,amsmath}
\usepackage{dsfont}
\usepackage{mathtools,amssymb}
\usepackage{extarrows}
\usepackage{tikz-cd}
\journal{\textnormal{}}
\usepackage{graphicx}
\usepackage{amsfonts}

\usepackage[colorlinks=true,linkcolor=blue,citecolor={blue}]{hyperref}

\theoremstyle{plain}
\newtheorem{theorem}{Theorem}[section]
\newtheorem{proposition}[theorem]{Proposition}
\newtheorem{lemma}[theorem]{Lemma}
\newtheorem{corollary}[theorem]{Corollary}
\newtheorem*{conjecture}{Conjecture}

\theoremstyle{definition}
\newtheorem{definition}{Definition}
\newtheorem{example}[theorem]{Example}
\newtheorem{remark}[theorem]{Remark}
\newtheorem*{Question}{Question}
\newtheorem*{Notations}{Notations}

\newcommand{\vecf}{\mathfrak{X}}
\newcommand{\met}{\langle\cdot\,,\cdot\rangle}
\newcommand{\rel}{\mathbf{R}}

\newcommand{\n}{\nabla}

\newcommand{\et}{\quad \textnormal{and}\quad}
\newcommand{\id}{\operatorname{Id}}

\newcommand{\tr}{\operatorname{tr}}
\newcommand{\A}{\mathbf{A}}

\newcommand{\dtz}{\dfrac{d}{dt}_{|_{t=0}}}

\newcommand{\g}{\mathfrak{g}}
\newcommand{\h}{\mathfrak{h}}
\newcommand{\m}{\mathfrak{m}}
\newcommand{\Ad}{\operatorname{Ad}}
\newcommand{\ad}{\operatorname{ad}}
\newcommand{\Span}{\operatorname{span}}

\newcommand{\T}{\operatorname{T}}

\newcommand{\Ll}{\operatorname{L}}

\newcommand{\ric}{\mathbf{ric}}
\newcommand{\rico}{\mathbf{Ric}}
\newcommand{\LS}{\operatorname{L}}

\numberwithin{equation}{section}

\newcommand{\somijn}{\displaystyle\sum_{i,j=1}^{n}}
\newcommand{\somijkn}{\displaystyle\sum_{i,j,k=1}^{n}}
\newcommand{\somi}{\displaystyle\sum_{i=1}^{n}}
\newcommand{\R}{\mathcal{R}}
\newcommand{\divr}{\delta\R}
\newcommand{\Rm}{\underline{\R}}
\newcommand{\divm}{\delta\Rm}
\newcommand{\cro}{[\,\cdot\,,\cdot\,]}
\newcommand{\p}{\mathfrak{p}}

\newcommand{\gl}{\mathfrak{gl}}

\newcommand{\rad}{\mathfrak{rad}}
\newcommand{\comment}[1]{}

\begin{document}
	
	\begin{frontmatter}
		
		\title{\textbf{Five-dimensional homogeneous Riemannian manifolds with harmonic curvature are Ricci-parallel}\vskip 0.5cm}
		
		\author[l2]{Charif Bourzik}
		\ead{bourzikcharif@gmail.com}
		\author[l1]{Othmane Dani}
		\ead{othmanedani@gmail.com}
		\author[l1]{Mehdi Nabil}
		\ead{mehdi1nabil@gmail.com}
		\address[l2]{Higher Institute of Nursing and Health Techniques Guelmim, Morocco}
		\address[l1]{Department of Mathematics, Faculty of Sciences and Technologies, Cadi Ayyad University  B.P.549 Gueliz Marrakesh, Morocco}
		
		\begin{abstract}
			
			In this paper, we show that $5$-dimensional homogeneous Riemannian manifolds with harmonic curvature are Ricci-parallel.
			
			\noindent{\bf Mathematics Subject Classification 2020:} 53B20; 53C07; 53C20; 53C25; 53C30. \\
			{\bf Keywords:} Homogeneous Riemannian manifolds, Harmonic curvature, Einstein-like metrics, Ricci-Codazzi, Ricci-parallel.
		\end{abstract}

	\end{frontmatter}	
	
	\section{Introduction and the statement of our main result}
	
	The study of curvature plays a central role in Riemannian geometry, as it provides a bridge between the local differential properties of a Riemannian manifold and its global geometric and topological structure (see, for example, \cite{Bess,Lee,Petr}). Among the numerous curvature conditions that have been studied (flatness, conformal flatness, local symmetry, the Einstein condition, Ricci parallelness, etc) the harmonicity of the curvature tensor occupies a distinguished position. Recall that a Riemannian manifold $(M,\met)$ is said to have \textit{harmonic} curvature if the divergence of its curvature tensor vanishes, i.e., $\delta\R=0$ (see \cite{Abr,Gry}). In the compact setting, this is equivalent to the fact that the Levi-Civita connection $\n$ is a critical point for the Yang-Mills functional in the space of all $\met$-metric connections (cf. \cite{Bess,Der1}). It is well known (e.g. \cite[p. 435]{Bess}) that a Riemannian manifold $(M,\met)$ has a harmonic curvature if and only if it is \textit{Ricci-Codazzi}, meaning that its Ricci tensor is a Codazzi tensor, i.e.,
	\begin{equation}\label{RicciCodazIntrod}
		\left(\n_X\ric\right)(Y,Z)=\left(\n_Y\ric\right)(X,Z),\qquad\forall\,X,Y,Z\in\vecf(M).
	\end{equation}
	Every Einstein Riemannian manifold (i.e., $\ric=\lambda\met$), every locally symmetric Riemannian manifold (i.e., $\n\R=0$), and more generally, every Riemannian manifold with parallel Ricci tensor (i.e., $\n\ric=0$) has harmonic curvature. Over the past decades, the harmonic curvature condition has been extensively investigated from geometric, analytic, and mathematical physics perspectives (see, for instance, \cite{Bess,Bourg',Bourg22,Datr,Gry,Li}). A number of rigidity results (e.g. \cite[p. 445]{Bess} and \cite[p. 280]{Bourg}) show that under suitable hypotheses such as compactness, positive sectional curvature, or nonzero signature in dimension four, harmonic curvature often forces the Riemannian manifold to be Einstein. Nevertheless, since harmonic curvature is too weak to imply either the Einstein condition or local symmetry, it is natural to ask whether it implies the parallelness of the Ricci tensor, and consequently (see \cite[p. 261]{Gry}), whether the Riemannian manifold is locally isometric to a Riemannian product of Einstein Riemannian manifolds. Unfortunately, many interesting compact examples (cf. \cite{Der0,Der}) demonstrate that harmonic curvature is strictly weaker than Ricci parallelness, making the class rich enough to exhibit diverse geometric phenomena.
	
	A Riemannian manifold $(M,\met)$ is called \textit{homogeneous} if the isometry Lie group $\operatorname{Iso}(M,\met)$ acts transitively on $M$. In \cite[p. 82]{Der1}, the following question was asked: 
	\begin{Question}
		Must a homogeneous Riemannian manifold with harmonic curvature
		have parallel Ricci tensor$?$
	\end{Question}
	Several affirmative answers have been obtained for various classes of homogeneous Riemannian manifolds, including nilpotent Riemannian Lie groups (cf. \cite{Datr}), solvable Riemannian Lie groups and Riemannian Lie groups of dimension at most six (cf. \cite{Abr}), certain compact Lie groups (cf. \cite{Wu}), homogeneous Riemannian manifolds of dimensions three and four (cf. \cite{Bu,Pod}), spheres and projective spaces (cf. \cite{Peng}), certain compact Riemannian symmetric spaces (cf. \cite{Qian}), certain generalized Wallach spaces (cf. \cite{LiCh}), naturally reductive homogeneous Riemannian manifolds (cf. \cite{Reb}), and conformally flat homogeneous Riemannian manifolds (cf. \cite{Tak}). Despite these advances, the question remains open in full generality. More recently, it was reformulated as the following conjecture (cf. \cite{Abr}):
	
	\begin{conjecture}
		Any homogeneous Riemannian manifold with harmonic curvature is Ricci-parallel.
	\end{conjecture}
	
	Since $5$-dimensional homogeneous Riemannian manifolds have attracted considerable attention in the literature (see, for example, \cite{Agrico,Ale,KowVan}), and since, to the best of our knowledge, nothing is known regarding this conjecture in this setting, the purpose of this paper is to provide an affirmative answer for this class of homogeneous Riemannian manifolds. More precisely, we prove the following theorem.
	
	\begin{theorem}\label{maintheo}
		Let $(M,\met)$ be a $5$-dimensional homogeneous Riemannian manifold with harmonic curvature. Then, its Ricci tensor must be parallel.
	\end{theorem}
	
	According to the discussion above (see also Corollary \ref{CorolApx}), to prove Theorem \ref{maintheo} it suffices to establish that:
	
	\begin{theorem}\label{maintheo2}
		Every $5$-dimensional Ricci-Codazzi homogeneous Riemannian manifold is Ricci-parallel.
	\end{theorem}
	
	Note that, in dimensions at most four, homogeneous Riemannian manifolds with harmonic curvature are locally symmetric (see \cite{Bu,Pod}). This phenomenon no longer holds in dimension five. For instance, the $5$-dimensional Stiefel manifold $\operatorname{SO}(4)/\operatorname{SO}(2)$ admits an Einstein $\operatorname{SO}(4)$-invariant metric (cf. \cite{Ale}), and hence has parallel Ricci tensor, although it is not locally symmetric (see \cite[p. 453]{Helg}).\vspace{3mm}
	
	\textit{Organization of the paper}. In Section \ref{Sect2}, we establish several algebraic properties of homogeneous Riemannian manifolds, culminating in Proposition \ref{propertiesreductive} and Corollary \ref{CorDim2}. In Section \ref{Sect3}, we provide algebraic characterizations of the Ricci-Codazzi condition (Theorem \ref{RicciCodazziTheorem}), the parallelness of the Ricci tensor (Theorem \ref{RicciParalcit1}), and the algebraic conditions under which a Ricci-Codazzi homogeneous Riemannian manifold has parallel Ricci tensor (Corollary \ref{CodazziParallel}). In Section \ref{Sect4}, we prove our main result, namely Theorem \ref{maintheo2}. Finally, the paper concludes with an appendix (Section \ref{Sect5}), where we present an elementary proof of the well-known equivalence between the harmonicity of the curvature tensor and the Ricci-Codazzi condition for Riemannian manifolds.
	
	\begin{Notations}
		We will always identify a homogeneous manifold $M$ with the quotient space $G/H$, where $G$ is a fixed Lie group acting transitively on $M$ and $H$ is a closed subgroup of $G$ (the isotropy subgroup at a fixed point of $M$). Moreover, we denote by $\g$ and $\h$ the Lie algebras of $G$ and $H$ respectively. Further, for a direct sum decomposition of a vector space $V=\mathfrak{a}\oplus\mathfrak{b}$, we represent the linear projection of an element $v\in V$ onto the vector subspace $\mathfrak{a}$ (resp. $\mathfrak{b}$) as $v_{\mathfrak{a}}$ (resp. $v_{\mathfrak{b}}$).
	\end{Notations}
	
	\section{Some facts about homogeneous Riemannian manifolds}\label{Sect2}
	
	In this section, we recall some basic facts about homogeneous Riemannian manifolds, many of which are well known in the literature. Let $(G/H,\met)$ be a homogeneous Riemannian $G$-space; that is, $G/H$ is a homogeneous $G$-space and $\met$ is a $G$-invariant Riemannian metric on it. Let $\g=\h\oplus\m$ be a fixed reductive decomposition (such a decomposition always exists; e.g. \cite[p. 210]{KowSz}), meaning that $\m$ is a vector subspace of $\g$ suh that $\Ad(H)(\m)=\m$. Following \cite[p. 200]{Kob2}, the Riemannian metric $\met$ gives rise to a positive definite inner product (also denoted by $\met$) on $\m$ which is $\Ad(H)$-invariant. On the other hand, according to \cite[p. 43]{Nom1}, since $G/H$ is a reductive homogeneous $G$-space, there is a one-to-one correspondence between the set of $G$-invariant connections on $G/H$ and the set of products on $\m$ which are $\Ad(H)$-invariant. Thus, since the Levi-Civita connection $\n$ of $(G/H,\met)$ is $G$-invariant, then under the identification of $T_H(G/H)$ with $\mathfrak{m}$, its associated product $\Ll:\mathfrak{m}\times\mathfrak{m}\rightarrow\mathfrak{m}$ is given by:
	\begin{equation}
		\Ll_x(y)=\A^\n_{x^\#}y^\#_H,
	\end{equation}    
	for $x,y\in\mathfrak{m}$, where $x^\#$ (resp. $y^\#$) is the fundamental vector field associated to $x$ (resp. $y$) whose complete flow is explicitly given by $(t,gH)\mapsto\exp_G(tx)gH$, and $\A^\n_{x^\#}:=\n_{x^\#}-\mathcal{L}_{x^\#}$ is the Nomizu's operator related to $\n$ and $x^\#$, where $\mathcal{L}$ denotes the Lie derivative. Hence, the Koszul's formula becomes:
	\begin{equation}\label{LeviCivita}
		\langle\Ll_x(y),z\rangle=\tfrac{1}{2}\Big\{\langle[x,y]_\m,z\rangle+\langle[z,x]_\m,y\rangle+\langle x,[z,y]_\m\rangle\Big\},
	\end{equation}
	for any $x,y,z\in\m$. Moreover, it is clear that the torsion $\mathcal{T}$, the curvature $\mathcal{R}$, and the Ricci $\ric$ tensor fields of $\n$ are also $G$-invariant. Thus, they are completely determined by their value at the origin $o:=H\in G/H$. Therefore, from \cite[pp. 45-47]{Nom1}, $\mathcal{T}$ and $\mathcal{R}$ can be expressed as follows:
	\begin{eqnarray}
		\mathcal{T}(x,y)&=&\Ll_x(y)-\Ll_y(x)-[x,y]_\m\label{Tor};\\
		\mathcal{R}(x,y)&=&\left[\Ll_x,\Ll_y\right]-\Ll_{[x,y]_\mathfrak{m}}-\ad_{[x,y]_\h}.\label{Curv}
	\end{eqnarray}
	Further, the Ricci tensor and the Ricci endomorphism are given by:
	\begin{equation*}
		\ric(x,y)=\somi \langle \R(x,e_i)y,e_i \rangle,\et \langle \rico(x),y\rangle=\ric(x,y),
	\end{equation*}
	for any $x,y\in\m$, where $\{e_i\}_{1\leqslant i\leqslant n:=\dim(G/H)}$ is an orthonormal basis of $(\m,\met)$. The fact that $\n$ is torsion-free and compatible with the metric $\met$, is equivalent respectively to:
	\begin{equation}\label{TorsionFreeMetric}
		\LS_x(y)-\LS_y(x)=[x,y]_\m,\et\langle\LS_x(y),z\rangle=-\langle y,\LS_x(z)\rangle.
	\end{equation}
	
	We begin by establishing some preliminary properties of reductive, not necessarily Riemannian, homogeneous manifolds.
	\begin{lemma}\label{LemInvVect}
		Let $G/H$ be a reductive homogeneous $G$-space with a fixed reductive decomposition $\g=\h\oplus\m$. Then, 
		\begin{itemize}
			\item[$1.$] The set $\m^\h:=\big\{x\in\m\,|\,\,[u,x]=0,\,\,\forall\,u\in\h\big\}$ consisting of $\h$-invariant vectors of $\m$ is stable under the product $\cro_\m$. Moreover, $(\m^\h,\cro_\m)$ is a Lie algebra;
			\item[$2.$] $\m^\h$ is isomorphic $($as Lie algebra$)$ to $\mathfrak{n}_\g(\h)/\h$, where $\mathfrak{n}_\g(\h)$ denotes the normalizer of $\h$ in $\g$, i.e., $\mathfrak{n}_\g(\h):=\big\{x\in\g\,|\,\,[u,x]\in\h,\,\,\forall\,u\in\h\big\}$;
			\item[$3.$] $[\m^\h,\m^\h]_\h\subseteq\mathfrak{z}(\h)$, where $\mathfrak{z}(\h)$ is the center of $\h$.
		\end{itemize}
	\end{lemma}
	
	\begin{proof}
		First, it is clear that $\m^\h$ is a vector subspace of $\m$. Moreover, for any $x,y\in\m^\h$ and $u\in\h$, we have:
		\begin{equation*}
			\big[u,[x,y]_\m\big]=\big[u,[x,y]\big]_\m=\big[[u,x],y\big]_\m+\big[x,[u,y]\big]_\m=0.
		\end{equation*}
		Therefore, $[x,y]_\m\in\m^\h$, and hence $\m^\h$ is a subalgebra of $(\m,\cro_\m)$. In addition, for $x,y,z\in\m^\h$, we compute:
		\begin{equation*}
			\big[x,[y,z]_\m\big]_\m=\big[x,[y,z]\big]_\m=\big[[x,y],z\big]_\m+\big[y,[x,z]\big]_\m=\big[[x,y]_\m,z\big]_\m+\big[y,[x,z]_\m\big]_\m.
		\end{equation*}
		This proves that $(\m^\h,\cro_\m)$ is a Lie algebra. On the other hand, if we consider the linear map $F:\mathfrak{n}_\g(\h)\to\m^\h,\,x\mapsto x_\m$, then one can easily check that $F$ is well-defined and surjective, since $\m^\h\subset\mathfrak{n}_\g(\h)$, furthermore $\ker(F)=\h$. We claim that $F$ is a Lie algebra homomorphism. Indeed, for $x,y\in\mathfrak{n}_\g(\h)$, we have:
		\begin{equation*}
			F([x,y])=[x,y]_\m=\underbrace{[x_\h,y]_\m}_{=0}+[x_\m,y]_\m=[x_\m,y_\m]_\m+\underbrace{[x_\m,y_\h]}_{=0}=\big[F(x),F(y)\big]_\m.
		\end{equation*}
		Thus, $\m^\h\cong\mathfrak{n}_\g(\h)/\h$. Finally, for $x,y\in\m^\h$ and $u\in\h$, we have:
		\begin{equation*}
			\big[[x,y]_\h,u\big]=\big[[x,y],u\big]=\big[[x,u],y\big]+\big[x,[y,u]\big]=0.
		\end{equation*}
		This completes the proof.
	\end{proof}
	
	The vector space $\m^\h$ plays a central role in the proof of our main theorem. The next two propositions describe the algebraic and topological consequences arising from the triviality or nontriviality of $\m^\h$ and are included for the sake of completeness. Recall that, the natural action of a Lie group $G$ on a homogeneous $G$-space $G/H$ is called almost
	effective if $H$ contains no connected normal subgroup of $G$ other than the identity element. It is well known that $G$ acts almost effectively on $G/H$ if and only if $\h$ contains no nonzero ideal of $\g$.
	
	\begin{proposition}
		Let $G/H$ be a homogeneous $G$-space where $G$ is a compact Lie group. If $G$ acts almost effectively on $G/H$ and $\h=\mathfrak{n}_\g(\h)$, then $\g$ is a semisimple Lie algebra. 
	\end{proposition}
	
	\begin{proof}
		First, the fact that $G$ is compact implies that $G/H$ is a reductive homogeneous $G$-space and so $\g=\h\oplus\m$ with $\Ad(H)(\m)=\m$, it also implies that the Lie algebra $\g$ is reductive, meaning that $\rad(\g)=\mathfrak{z}(\g)$. Therefore, to show that $\g$ is semisimple, it suffices to establish that $\mathfrak{z}(\g)=\{0\}$. For any $z\in\mathfrak{z}(\g)$ and $u\in\h$, we have:
		\begin{equation*}
			0=[z,u]=\underbrace{[z_\h,u]}_{\in\h}-\underbrace{[u,z_\m]}_{\in\m}.
		\end{equation*}
		Therefore, $z_\h\in\mathfrak{z}(\h)$ and $z_\m\in\m^\h$. However, from Lemma \ref{LemInvVect}, we have $\m^\h=\{0\}$. Hence, $\mathfrak{z}(\g)\subseteq\h$, which gives that $\mathfrak{z}(\g)$ must be equal to zero, since $G$ acts almost effectively on $G/H$. 
	\end{proof}
	
	\begin{proposition}\label{Euler=0}
		Let $G/H$ be a compact reductive homogeneous $G$-space with connected isotropy Lie subgroup $H$. If $\h$ is strictly contained in $\mathfrak{n}_\g(\h)$, then $G/H$ has zero Euler characteristic, i.e. $\chi(G/H)=0$.
	\end{proposition}
	
	\begin{proof}
		Let $\g=\h\oplus\m$ be a reductive decomposition of $\g$. Since $H$ is connected, we have:
		\begin{equation*}
			\m^\h=\m^H:=\big\{x\in\m\,|\,\,\Ad_h(x)=x,\,\,\forall\,h\in H\big\}.
		\end{equation*}
		In addition, using the fact that $\h$ is strictly contained in $\mathfrak{n}_\g(\h)$, one can deduce from Lemma \ref{LemInvVect} that $\m^\h\neq\{0\}$. Take $x\in\m^\h\backslash\{0\}$ and define a non-vanishing vector field $X\in\vecf(G/H)$ as follows:
		\begin{equation*}
			X_{gH}:=\dtz g\exp_G(tx)H,
		\end{equation*}
		for all $gH\in G/H$. Thus, the result follows from the Poincaré-Hopf theorem.
	\end{proof}
	
	\begin{example}[\textbf{Stiefel manifolds}]
		Given $k,n\in\mathbb{N}$ such that $2\leqslant k\leqslant n$, let $\mathcal{S}_{k,n}$ be the set of all orthonormal $k$-frames in the Euclidean space $(\rel^n,\met_0)$. It is well known that $\mathcal{S}_{k,n}$ is diffeomorphic to the homogeneous manifold $\operatorname{SO}(n)/\operatorname{SO}(n-k)$. Moreover, $\operatorname{SO}(n)/\operatorname{SO}(n-k)$ is reductive, because $\mathfrak{so}(n)=\mathfrak{so}(n-k)\oplus\m$, where:
		\begin{equation*}
			\mathfrak{so}(n-k)\simeq\left\{\begin{pmatrix}
				\textbf{0}&\textbf{0}\\\textbf{0}&S
			\end{pmatrix}\,\big|\,\,S\in\mathfrak{so}(n-k)\right\},$$
			$$\m=\left\{\begin{pmatrix}
				X&-Y^{\operatorname{T}}\\Y&\textbf{0}
			\end{pmatrix}\,\big|\,\,X\in\mathfrak{so}(n-k),\,Y\in\mathcal{M}_{n-k,k}(\rel)\right\}.
		\end{equation*}
		Since $\m^{\operatorname{SO}(n-k)}=\Bigl\{\begin{psmallmatrix}
			X&\textbf{0}\\\textbf{0}&\textbf{0}
		\end{psmallmatrix}\,\big|\,\,X\in\mathfrak{so}(n-k)\Bigr\}\neq\{\textbf{0}\}$, then $\chi\big(\mathcal{S}_{k,n}\big)=0$ according to Proposition \ref{Euler=0}.
	\end{example}
	
	The next lemma provides a key decomposition of the vector space $\m$, which will be used throughout the remainder of the paper.
	
	\begin{lemma}\label{lemDecompM}
		Let $G/H$ be a reductive homogeneous $G$-space with completely reducible isotropy linear representation $\ad^{\g/\h}:\h\to\gl(\g/\h)$. For any reductive decomposition $\g=\h\oplus\m$ of $\g$, we have:
		$$\m=\m^\h\oplus[\h,\m].$$ In addition, the following inclusion holds $\big[\m^\h,[\h,\m]\big]_\m\subseteq [\h,\m]$.
	\end{lemma}
	
	\begin{proof}
		Let $\g=\h\oplus\m$ be a fixed reductive decomposition. Since $\ad^\m:\h\to\gl(\m)$ is completely reducible and $\m^\h$ is an $\h$-submodule of $\m$, there exists an $\h$-submodule complement $\m_1$, so that $\m=\m^\h\oplus\m_1$. Using the fact that $\m^\h$ is a trivial $\h$-module, one can easily see that $[\h,\m]\subseteq\m_1$. In particular, $\m^\h\cap[\h,\m]=\{0\}$. We claim that $\m_1=[\h,\m]$. In fact, since $\m^\h\oplus[\h,\m]$ is an $\h$-submodule of $\m$, one can write $\m=\m^\h\oplus[\h,\m]\oplus\m_0$, where $\m_0$ is an $\h$-submodule of $\m$. Moreover, we have:
		\begin{equation*}
			[\h,\m_0]\subseteq [\h,\m]\cap\m_0=\{0\}.
		\end{equation*}
		This shows that $\m_0\subseteq\m^\h$, and therefore $\m_0=\{0\}$. Finally, for any $u\in\h,\,x\in\m^\h,\,$ and $y\in\m$, we have:
		\begin{equation*}
			\big[x,[u,y]\big]_\m=\big[u,[x,y]\big]_\m=\big[u,[x,y]_\m\big].
		\end{equation*}
		As a result, we get the desired inclusion $\big[\m^\h,[\h,\m]\big]_\m\subseteq [\h,\m]$.
	\end{proof}
	
	Applying this to our context of homogeneous Riemannian manifolds, we obtain the following proposition.
	
	\begin{proposition}\label{propertiesreductive}
		Let $(G/H,\met)$ be a homogeneous Riemannian $G$-space with a fixed reductive decomposition $\g=\h\oplus\m$. The following assertions hold:
		\begin{enumerate}
			\item[$1.$] The isotropy linear representation $\ad^\m:\h\to\gl(\m)$ is completely reducible, and $\m=\m^\h\overset{\perp}{\oplus}[\h,\m]$;
			\item[$2.$] Any $1$-dimensional $\h$-submodule of $\m$ is contained in $\m^\h$. In particular, if $\p$ is a nonzero $\h$-submodule of $[\h,\m]$, then $\dim\p\geqslant 2$;
			\item[$3.$] Any $\h$-submodule of $[\h,\m]$ with dimension less than $4$ is irreducible;
			\item[$4.$] If $G$ acts almost effectively on $G/H$, then $\h$ is a reductive Lie algebra, meaning that $\rad(\h)=\mathfrak{z}(\h)$.
		\end{enumerate}
	\end{proposition}
	
	\begin{proof}
		First, to show that the isotropy linear representation $\ad^\m:\h\to\gl(\m)$ is completely reducible, it suffices to prove that each $\h$-submodule of $\m$ possesses an $\h$-submodule complement. Let $\p\subseteq\m$ be an $\h$-submodule of $\m$ and denote by $\p^{\perp_{\met}}\subseteq\m$ its orthogonal complement with respect to $\met$. Clearly that $\m=\p\oplus\p^{\perp_{\met}}$, and for any $x\in\p^{\perp_{\met}},\,y\in\p,\,u\in\h$ we have:
		\begin{equation*}
			\langle[u,x],y\rangle=-\langle x,[u,y]\rangle=0.
		\end{equation*} 
		Hence, $\p^{\perp_{\met}}$ is an $\h$-submodule of $\m$. The decomposition $\m=\m^\h\oplus[\h,\m]$ follows from Lemma \ref{lemDecompM}, and one can easily verify, since $\met$ is $\h$-invariant, that $[\h,\m]^{\perp_{\met}}=\m^\h$. For the second assertion, let $\p$ be a $1$-dimensional $\h$-submodule of $\m$. Write $\p=\rel x_0$ for some $x_0\in\m\backslash\{0\}$ with $\langle x_0,x_0\rangle=1$. Since $\p$ is a $1$-dimensional $\h$-module, there exists a linear form $\alpha\in\h^*$ such that:
		\begin{equation*}
			[u,x_0]=\alpha(u)x_0\qquad\forall\,u\in\h.
		\end{equation*}
		On the other hand, using the fact that $\met$ is $\h$-invariant, we compute:
		\begin{equation*}
			2\alpha(u)=2\alpha(u)\langle x_0,x_0\rangle=\langle[u,x_0],x_0\rangle+\langle x_0,[u,x_0]\rangle=0,
		\end{equation*}
		for all $u\in\h$. Therefore, $\alpha=0$, and then $[\h,\p]=\{0\}$. Thus, $\p\subseteq\m^\h$. The third assertion follows from the second one. In fact, if $\p$ is a two or a three dimensional $\h$-submodule of $[\h,\m]$ which is not irreducible, then $\p$ must contain a $1$-dimensional $\h$-submodule of $\m$. However, this is impossible, since $\p\subseteq [\h,\m]$. Finally, using the fact that the isotropy linear representation $\ad^\m:\h\to\gl(\m)$ is faithful and completely reducible, one can deduce from \cite[p. 56]{Bourb} that $\h$ is reductive.
	\end{proof}
	
	We end this section by the following corollary.
	
	\begin{corollary}\label{CorDim2}
		Let $(G/H,\met)$ be a homogeneous Riemannian $G$-space with a fixed reductive decomposition $\g=\h\oplus\m$. Then, for any $2$-dimensional $\h$-submodule $\p$ of $\m$, and any orthonormal basis $\{x,y\}$ of $\p$, there exists a linear form $\alpha\in\h^*$ such that:
		\begin{equation}\label{OrthoBasisDim2}
			[u,x]=\alpha(u)y,\et [u,y]=-\alpha(u)x,\qquad\forall\,u\in\h.
		\end{equation}
		Moreover, $[\p,\p]\subseteq\mathfrak{z}(\h)\oplus\m^\h$.
	\end{corollary}
	
	\begin{proof}
		The first assertion follows from the fact that $\ad^\m(\h)\subseteq\mathfrak{so}(\p,\met)$ and $\p$ is $2$-dimensional. To prove the inclusion, it suffices to show that $[x,y]\in\mathfrak{z}(\h)\oplus\m^\h$, where $\{x,y\}$ is an orthonormal basis of $\p$. For any $u\in\h$, we compute:
		\begin{equation*}
			\big[u,[x,y]\big]=\big[[u,x],y\big]+\big[x,[u,y]\big]\stackrel{\rm\eqref{OrthoBasisDim2}}{=}0.
		\end{equation*}
		Thus $\big[u,[x,y]_\m\big]=\big[u,[x,y]_\h\big]=0$, consequently $[x,y]_\m\in\m^\h$ and $[x,y]_\h\in\mathfrak{z}(\h)$.
	\end{proof}

	\section{Ricci-Codazzi homogeneous Riemannian manifolds}\label{Sect3}
	
	The aim of this section is to establish an algebraic characterization of the parallelness of the Ricci tensor for homogeneous Riemannian manifolds. In what follows, $(G/H,\met)$ will be a homogeneous Riemannian $G$-space, $\g=\h\oplus\m$ a fixed reductive decomposition, and $\m=\m^\h\oplus[\h,\m]$ the orthogonal decomposition of $\m$ introduced in Proposition \ref{propertiesreductive}. Since, $\rico:\m\to\m$ is a symmetric linear endomorphism of $\m$ and $[\h,\m]^{\perp_{\met}}=\m^\h$, using the $\h$-invariance of $\ric$, we have:
	\begin{equation}\label{RicStable}
		\rico(\m^\h)\subseteq\m^\h,\quad\et\quad \rico([\h,\m])\subseteq [\h,\m].
	\end{equation}
	In particular, we denote by $\rico_{[\h,\m]}:[\h,\m]\to[\h,\m]$ the restriction of $\rico:\m\to\m$ to $[\h,\m]$. Further, using the symmetry of $\rico:\m\to\m$, the spectral
	theorem allows us to write an orthogonal direct sum decomposition
	\begin{equation}\label{spect}
		\m=\m_1\overset{\perp}{\oplus}\cdots\overset{\perp}{\oplus}\m_r,
	\end{equation}
	where $r\geqslant 1$, $\m_i:=\ker\big(\rico-\lambda_i\id_\m\big)$ for $i\in\{1,\ldots,r\}$, and $\lambda_1,\ldots,\lambda_r$ are the pairwise distinct real eigenvalues of $\rico:\m\to\m$. The decomposition \eqref{spect} is called the spectral decomposition of $\rico:\m\to\m$. In addition, it is clear that all $\m_i$'s are $\h$-submodules of $\m$. The Ricci-Codazzi condition \eqref{RicciCodazIntrod} is equivalent to:
	\begin{equation}\label{CodazziCond0}
		\big(\Ll_{x}\ric\big)(y,z)=\big(\Ll_{y}\ric\big)(x,z),
	\end{equation}
	for any $x,y,z\in\m$, where $\big(\Ll_{x}\ric\big)(y,z):=\ric\big(\Ll_x(y),z\big)+\ric\big(y,\Ll_x(z)\big)$.
	In terms of eigenspaces $\m_i$'s, \eqref{CodazziCond0} can be written as:
	\begin{equation*}
		\big(\Ll_{x_i}\ric\big)(y_j,z_k)=	\big(\Ll_{y_j}\ric\big)(x_i,z_k),\qquad\forall\,x_i\in\m_i,y_j\in\m_j,z_k\in\m_k.
	\end{equation*}
	
	Another remarkable characterization of the Ricci-Codazzi condition is given by the following lemma.
	
	\begin{lemma}
		$(G/H,\met)$ is Ricci-Codazzi if and only if 
		\begin{equation}\label{CodazziCond}
			\rico\big([x,y]_\m\big)=\Ll_x\big(\rico(y)\big)-\Ll_y\big(\rico(x)\big),
		\end{equation}
		for any $x,y\in\m$.
	\end{lemma}
	
	\begin{proof}
		We have that $(G/H,\met)$ is Ricci-Codazzi if and only if \eqref{CodazziCond0} holds. In other words,
		\begin{equation*}
			\ric\big(\Ll_x(y),z\big)+\ric\big(y,\Ll_x(z)\big)=\ric\big(\Ll_y(x),z\big)+\ric\big(x,\Ll_y(z)\big).
		\end{equation*}
		Equivalently, 
		\begin{equation*}
			\ric\big(\Ll_x(y),z\big)-\ric\big(\Ll_y(x),z\big)=\ric\big(x,\Ll_y(z)\big)-\ric\big(y,\Ll_x(z)\big).
		\end{equation*}
		Rewriting this, using the fact that $\n$ is torsion-free, we get
		\begin{equation*}
			\langle\rico\big([x,y]_\m\big),z\rangle=\langle\rico(x),\Ll_y(z)\rangle-\langle\rico(y),\Ll_x(z)\rangle.
		\end{equation*}
		Hence, the equivalence follows using \eqref{TorsionFreeMetric}.
	\end{proof}
	
	Consequently, we obtain the next corollary.
	
	\begin{corollary}\label{eigenspacesstable}
		If $(G/H,\met)$ is Ricci-Codazzi, then for any $i\in\{1,\ldots,r\}$, the following inclusions hold:
		$$[\m_i,\m_i]_\m\subseteq\m_i,\et\Ll(\m_i,\m_i)\subseteq\m_i.$$
	\end{corollary}
	
	\begin{proof}
		The first one follows directly from \eqref{CodazziCond} and \eqref{TorsionFreeMetric}. The second one follows from \eqref{LeviCivita} and \eqref{CodazIeqlJ} (see Remark \ref{RemkCoroProf} below).
	\end{proof}
	
	Before establishing the main algebraic characterization of Ricci-Codazzi homogeneous Riemannian manifolds, we first recall that every Einstein homogeneous Riemannian manifold has parallel Ricci tensor.
	
	\begin{proposition}\label{EinstParallel}
		If $\rico:\m\to\m$ has only one real eigenvalue, i.e., $\rico=\lambda\id_\m$, then $(G/H,\met)$ is Ricci-parallel.
	\end{proposition}
	
	\begin{proof}
		The parallelness of the Ricci tensor follows from the parallelness of the Riemannian metric $\met$. 
	\end{proof}
	
	The following theorem provides the main algebraic characterization of Ricci-Codazzi homogeneous Riemannian manifolds (cf. \cite[p. 5]{Reb}).
	
	\begin{theorem}\label{RicciCodazziTheorem}
		Let $(G/H,\met)$ be a homogeneous Riemannian $G$-space, and fix a reductive decomposition $\g=\h\oplus\m$. Let $\m=\m_1\oplus\cdots\oplus\m_r$ be the spectral decomposition of $\rico:\m\to\m$. Then, $(G/H,\met)$ is Ricci-Codazzi if and only if 
		\begin{equation}\label{EquivCodazzi}
			(\lambda_i-\lambda_k)^2\langle [x_i,y_j]_\m,z_k\rangle+(\lambda_i-\lambda_j)^2\langle[x_i,z_k]_\m,y_j\rangle=0,
		\end{equation}
		for all $x_i\in\m_i,\,y_j\in\m_j,$ and $z_k\in\m_k$, where $i,j,k\in\{1,\ldots,r\}$.
	\end{theorem}
	
	\begin{proof}
		Assume that $(G/H,\met)$ is Ricci-Codazzi. Taking $x_i\in\m_i,\,y_j\in\m_j$, and $z_k\in\m_k$, we compute: 
		\begin{eqnarray*}
			2\lambda_k\langle[x_i,y_j]_\m,z_k\rangle
			&\stackrel{\rm\eqref{CodazziCond}}{=}&2\lambda_i\langle\Ll_{y_j}(z_k),x_i\rangle-2\lambda_j\langle\Ll_{x_i}(z_k),y_j\rangle\\
			&\stackrel{\rm\eqref{LeviCivita}}{=}&(\lambda_i+\lambda_j)\langle[x_i,y_j]_\m,z_k\rangle+(\lambda_i-\lambda_j)\langle x_i,[y_j,z_k]_\m\rangle\\
			&&+(\lambda_i-\lambda_j)\langle y_j,[x_i,z_k]_\m\rangle.
		\end{eqnarray*}
		In other words, 
		\begin{equation}\label{Codproof}
			(2\lambda_k-\lambda_i-\lambda_j)\langle[x_i,y_j]_\m,z_k\rangle-(\lambda_i-\lambda_j)\Big(\langle x_i,[y_j,z_k]_\m\rangle+\langle y_j,[x_i,z_k]_\m\rangle\Big)=0.\quad
		\end{equation}
		Permuting $(i,j,k)\mapsto (k,i,j)$ and $(i,j,k)\mapsto (j,k,i)$ in \eqref{Codproof}, we get:
		\begin{equation}\label{Codproof1}
			\left\{
			\begin{array}{r c l}
				(2\lambda_j-\lambda_k-\lambda_i)\langle[z_k,x_i]_\m,y_j\rangle-(\lambda_k-\lambda_i)\Big(\langle z_k,[x_i,y_j]_\m\rangle+\langle x_i,[z_k,y_j]_\m\rangle\Big)&=&0,\\\\ 
				(2\lambda_i-\lambda_j-\lambda_k)\langle[y_j,z_k]_\m,x_i\rangle-(\lambda_j-\lambda_k)\Big(\langle y_j,[z_k,x_i]_\m\rangle+\langle z_k,[y_j,x_i]_\m\rangle\Big)&=&0.\hfill
			\end{array}
			\right.
		\end{equation}
		Thus, if we put:
		\begin{equation*}
			\mu:=\lambda_i-\lambda_j,\,\,\nu:=\lambda_j-\lambda_k,\,\,a:=\langle x_i,[y_j,z_k]_\m\rangle,\,\,b:=\langle z_k,[x_i,y_j]_\m\rangle,\,\,\textnormal{and}\,\,c:=\langle y_j,[z_k,x_i]_\m\rangle.
		\end{equation*}
		Then, rearranging terms in \eqref{Codproof} and \eqref{Codproof1}, we obtain:
		\begin{equation*}
			\begin{pmatrix}
				\mu& 2\nu+\mu& -\mu\\
				\mu+\nu&-(\mu+\nu)&\nu-\mu\\
				2\mu+\nu&\nu&-\nu
			\end{pmatrix}\begin{pmatrix}
				a\\
				b\\
				c
			\end{pmatrix}=\begin{pmatrix}
				0\\ 0\\ 0
			\end{pmatrix}.
		\end{equation*} 
		If $i=j$, then the equality \eqref{Codproof} is equivalent to \eqref{EquivCodazzi}. If $i\neq j$, then $\mu\neq 0$, and therefore, by Gaussian elimination, we get:
		\begin{equation*}
			\begin{pmatrix}
				\mu& 2\nu+\mu& -\mu\\
				0&-\tfrac{2(\mu+\nu)^2}{\mu}&2\nu\\
				0&0&2(\mu-\nu)
			\end{pmatrix}\begin{pmatrix}
				a\\
				b\\
				c
			\end{pmatrix}=\begin{pmatrix}
				0\\ 0\\ 0
			\end{pmatrix}.
		\end{equation*} 
		Hence, $\nu c=\mu c$ and $\mu\nu c=(\mu+\nu)^2b$. As a result,
		\begin{eqnarray*}
			(\lambda_i-\lambda_k)^2\langle [x_i,y_j]_\m,z_k\rangle+(\lambda_i-\lambda_j)^2\langle[x_i,z_k]_\m,y_j\rangle&=&(\mu+\nu)^2b-\mu^2c\\
			&=&(\mu+\nu)^2b-\mu\nu c\\
			&=&0.
		\end{eqnarray*}
		Conversely, assume that the equality \eqref{EquivCodazzi} holds, and let us show that $(G/H,\met)$ is Ricci-Codazzi, which is equivalent to establishing \eqref{Codproof}. For any  $x_i\in\m_i,\,y_j\in\m_j,$ and $z_k\in\m_k$, we have:
		\begin{equation*}
			\left\{
			\begin{array}{r c l}
				(\lambda_i-\lambda_k)^2\langle [x_i,y_j]_\m,z_k\rangle+(\lambda_i-\lambda_j)^2\langle[x_i,z_k]_\m,y_j\rangle&=&0,\\\\ 
				(\lambda_j-\lambda_i)^2\langle [y_j,z_k]_\m,x_i\rangle+(\lambda_j-\lambda_k)^2\langle[y_j,x_i]_\m,z_k\rangle&=&0.\hfill
			\end{array}
			\right.
		\end{equation*}
		Thus, by summing the two equalities above, we get:
		\begin{equation*}
			(\lambda_j-\lambda_i)(2\lambda_k-\lambda_i-\lambda_j)\langle[x_i,y_j]_\m,z_k\rangle+(\lambda_j-\lambda_i)^2\Big(\langle x_i,[y_j,z_k]_\m\rangle+\langle y_j,[x_i,z_k]_\m\rangle\Big)=0.\quad
		\end{equation*}
		Since when $i=j$, the equality \eqref{EquivCodazzi} is equivalent to \eqref{Codproof}, it suffices to consider the case $i\neq j$. However, after simplifying by $(\lambda_j-\lambda_i)$ in the equality above, we obtain \eqref{Codproof}. This completes the proof.
	\end{proof}
	
	\begin{remark}\label{RemkCoroProf}
		In the case where $k=j$, we deduce from \eqref{EquivCodazzi} that: 
		\begin{equation}\label{CodazIeqlJ}
			\langle[x_i,y_j]_\m,z_j\rangle+\langle[x_i,z_j]_\m,y_j\rangle=0
		\end{equation}
	\end{remark}
	
	The following theorem was established in \cite[p. 19]{DerSh} (see also \cite[p. 438]{Bess}).
	
	\begin{theorem}
		Let $(G/H,\met)$ be a Ricci-Codazzi homogeneous Riemannian $G$-space, $\g=\h\oplus\m$ a fixed reductive decomposition, and $\m=\m_1\oplus\cdots\oplus\m_r$ the spectral decomposition of $\rico:\m\to\m$. Then, for any $x_i\in\m_i,\,y_j\in\m_j,\,z_k\in\m_k$ such that $k\notin\{i,j\}$, we have:
		\begin{equation}\label{Curv=0}
			\R(x_i,y_j)z_k=0.
		\end{equation}
	\end{theorem}
	
	The next result is an extension of Proposition \ref{EinstParallel}.
	
	\begin{proposition}\label{Twoeigenval}
		If $(G/H,\met)$ is Ricci-Codazzi and $\rico:\m\to\m$ has exactly two eigenvalues, then $(G/H,\met)$ is Ricci-parallel.
	\end{proposition}
	
	\begin{proof}
		Assume that $\rico:\m\to\m$ has exactly two eigenvalues $\lambda_1,\lambda_2\in\rel$, and let $\m=\m_1\oplus\m_2$ be the spectral decomposition of $\rico:\m\to\m$. For $x_i\in\m_i,\,y_j\in\m_j,$ and $z_k\in\m_k$ with $i,j,k\in\{1,2\}$, we compute:
		\begin{eqnarray*}
			\big(\Ll_{x_i}\ric\big)(y_j,z_k)&=&\ric\big(\Ll_{x_i}(y_j),z_k\big)+\ric\big(y_j,\Ll_{x_i}(z_k)\big)\\
			&=&\lambda_k\langle\Ll_{x_i}(y_j),z_k\rangle+\lambda_j\langle y_j,\Ll_{x_i}(z_k)\rangle\\
			&=&(\lambda_k-\lambda_j)\langle\Ll_{x_i}(y_j),z_k\rangle\\
			&=&0,
		\end{eqnarray*}
		where the last equality follows from the fact that $\Ll(\m_i,\m_i)\subseteq\m_i$.
	\end{proof}
	
	We now turn our attention to the algebraic characterization of the parallelness of the Ricci tensor. Recall that, $(G/H,\met)$ is Ricci-parallel if and only if $\big(\Ll_{x}\ric\big)(y,z)=0$ for all $x,y,z\in\m$.
	
	\begin{proposition}\label{propRicciParal}
		$(G/H,\met)$ is Ricci-parallel if and only if $\Ll(\m_i,\m_j)\subseteq\m_j$ for any $i,j\in\{1,\ldots,r\}$.
	\end{proposition}
	
	\begin{proof}
		For $x_i\in\m_i,\,y_j\in\m_j,\,z_k\in\m_k$ with $i,j,k\in\{1,\ldots,r\}$ and $k\neq j$, we have:
		\begin{eqnarray*}
			\big(\Ll_{x_i}\ric\big)(y_j,z_k)&=&\ric\big(\Ll_{x_i}(y_j),z_k\big)+\ric\big(y_j,\Ll_{x_i}(z_k)\big)\\
			&=&\lambda_k\langle\Ll_{x_i}(y_j),z_k\rangle+\lambda_j\langle y_j,\Ll_{x_i}(z_k)\rangle\\
			&=&(\lambda_k-\lambda_j)\langle\Ll_{x_i}(y_j),z_k\rangle.
		\end{eqnarray*}
		Hence, the equivalence follows.
	\end{proof}
	
	Another noteworthy characterization of the Ricci parallelness is:
	
	\begin{theorem}\label{RicciParalcit1}
		Let $(G/H,\met)$ be a homogeneous Riemannian $G$-space and fix a reductive decomposition $\g=\h\oplus\m$. Let $\m=\m_1\oplus\cdots\oplus\m_r$ be the spectral decomposition of $\rico:\m\to\m$. Then, $(G/H,\met)$ is Ricci-parallel if and only if $[\m_i,\m_i]_\m\subseteq\m_i$, $[\m_i,\m_j]_\m\subseteq\m_i\oplus\m_j$, and
		\begin{equation}\label{CondRicParal}
			\langle[x_i,y_j]_\m,z_j\rangle+\langle[x_i,z_j]_\m,y_j\rangle=0,
		\end{equation}
		for all $x_i\in\m_i,\,y_j,\,z_j\in\m_j,$ with $i,j\in\{1,\ldots,r\}$ such that $i\neq j$.
	\end{theorem}
	
	\begin{proof}
		Suppose that $(G/H,\met)$ is Ricci-parallel, then using \eqref{TorsionFreeMetric} and Proposition \ref{propRicciParal}, we have: 
		\begin{equation*}
			[\m_i,\m_j]_\m\subseteq \Ll(\m_j,\m_i)+\Ll(\m_i,\m_j)\subseteq\m_i+\m_j,
		\end{equation*}
		for all $i,j\in\{1,\ldots,r\}$. Moreover, for $x_i\in\m_i$ and $y_j,\,z_j\in\m_j$ with $i\neq j$, we have:
		\begin{eqnarray*}
			2\langle[x_i,y_j]_\m,z_j\rangle&=&2\langle\Ll_{x_i}(y_j),z_j\rangle-\underbrace{\langle\Ll_{y_j}(x_i),z_j\rangle}_{=0}\\
			&\stackrel{\rm\eqref{LeviCivita}}{=}&\langle[x_i,y_j]_\m,z_j\rangle+\langle[z_j,x_i]_\m,y_j\rangle+\underbrace{\langle x_i,[z_j,y_j]_\m\rangle}_{=0}\\
			&=&\langle[x_i,y_j]_\m,z_j\rangle-\langle [x_i,z_j]_\m,y_j\rangle.
		\end{eqnarray*}
		Thus, \eqref{CondRicParal} follows. Conversely, assume that the conditions of the theorem hold, and let us prove that $(G/H,\met)$ is Ricci-parallel. For $x_i,\,y_i\in\m_i,\,z_k\in\m_k$ with $i,k\in\{1,\ldots,r\}$ and $k\neq i$, we have:
		\begin{equation*}
			\langle\Ll_{x_i}(y_i),z_k\rangle=\tfrac{1}{2}\Big\{\underbrace{\langle[x_i,y_i]_\m,z_k\rangle}_{=0}+\langle[z_k,x_i]_\m,y_i\rangle+\langle [z_k,y_i]_\m,x_i\rangle\Big\}\stackrel{\rm\eqref{CondRicParal}}{=}0.
		\end{equation*}
		Therefore, $\Ll(\m_i,\m_i)\subseteq\m_i$. Similarly, for $x_i\in\m_i,\,y_j\in\m_j,\,z_k\in\m_k$ with $i\neq j$ and $k\notin\{i,j\}$, we have: 
		\begin{equation*}
			\langle\Ll_{x_i}(y_j),z_k\rangle=\tfrac{1}{2}\Big\{\underbrace{\langle[x_i,y_j]_\m,z_k\rangle}_{=0}+\underbrace{\langle[z_k,x_i]_\m,y_j\rangle}_{=0}+\underbrace{\langle [z_k,y_j]_\m,x_i\rangle}_{=0}\Big\}=0.
		\end{equation*}
		Hence, $\Ll(\m_i,\m_j)\subseteq\m_i
		\oplus\m_j$. In addition, for $k=i$, we have:
		\begin{equation*}
			\langle\Ll_{x_i}(y_j),z_i\rangle=\tfrac{-1}{2}\Big\{\langle[y_j,x_i]_\m,z_i\rangle+\underbrace{\langle[x_i,z_i]_\m,y_j\rangle}_{=0}+\langle [y_j,z_i]_\m,x_i\rangle\Big\}\stackrel{\rm\eqref{CondRicParal}}{=}0.
		\end{equation*}
		This shows that $\Ll(\m_i,\m_j)\subseteq\m_j$, which finishes the proof.
	\end{proof}
	
	Applying this to Ricci-Codazzi homogeneous Riemannian manifolds, we deduce that:
	
	\begin{corollary}\label{CodazziParallel}
		Let $(G/H,\met)$ be a Ricci-Codazzi homogeneous Riemannian $G$-space, $\g=\h\oplus\m$ a fixed reductive decomposition, and $\m=\m_1\oplus\cdots\oplus\m_r$ the spectral decomposition of $\rico:\m\to\m$. Then, $(G/H,\met)$ is Ricci-parallel if and only if $[\m_i,\m_j]_\m\subseteq\m_i\oplus\m_j$.
	\end{corollary}
	
	\begin{proof}
		The equivalence follows from Corollary \ref{eigenspacesstable}, Remark \ref{RemkCoroProf}, and Theorem \ref{RicciParalcit1}.  
	\end{proof}
	
	\section{The proof of Theorem \ref{maintheo2}}\label{Sect4}
	
	Throughout this section, let $(G/H,\met)$ be a $5$-dimensional Ricci-Codazzi homogeneous Riemannian $G$-space, and let $\g=\h\oplus\m$ be a fixed reductive decomposition. Using the notation introduced in the previous sections, we prove our main result by considering the different possible dimensions of the vector space $\m^\h$. The proof is based primarily on Corollary \ref{CodazziParallel}. Since $G/H$ is $5$-dimensional and $\m^\h$ is a vector subspace of $\m$, we have $\dim\m^\h\in\{0,1,2,3,4,5\}$. If $\dim [\h,\m]=0$, then $\dim\m^\h=5$, and $(\m,\cro_\m,\met)$ is a $5$-dimensional Riemannian Lie algebra (see Lemma \ref{LemInvVect}) which satisfies the Codazzi condition, namely \eqref{EquivCodazzi}. Thus, the result follows from \cite[p. 8]{Abr}. On the other hand, if $\dim [\h,\m]\neq 0$, then Proposition \ref{propertiesreductive} implies that $\dim\m^\h\leqslant 3$. Therefore, it remains to consider the following four cases: $\dim\m^\h\in\{0,1,2,3\}$.

	\subsection{The case where $\m^\h=\{0\}$}
	
	In this case $[\h,\m]$ is $5$-dimensional, and hence by  Proposition \ref{propertiesreductive}, it is either an irreducible $\h$-submodule of $\m$ or splits into two irreducible $\h$-submodules of $\m$ one is $2$-dimensional and the other is $3$-dimensional. As consequence, $\rico:\m\to\m$ has at most two eigenvalues. Therefore, by Proposition \ref{EinstParallel} and Proposition \ref{Twoeigenval}, we conclude that $\n\ric=0$.
	
	\subsection{The case where $\m^\h$ is $1$-dimensional}
	
	In this case $[\h,\m]$ is $4$-dimensional. Therefore, by Proposition \ref{propertiesreductive}, it is either an irreducible $\h$-submodule of $\m$ or splits into two $2$-dimensional irreducible $\h$-submodules of $\m$. In the first case, $\rico:\m\to\m$ will have at most two eigenvalues, and hence the result follows from Propositions \ref{EinstParallel} and \ref{Twoeigenval}. In the second case, $\rico:\m\to\m$ will have at most three eigenvalues. Thus, the only remaining case to consider is when $\rico:\m\to\m$ has exactly three distinct eigenvalues. Assume that $\rico:\m\to\m$ has three mutually distinct eigenvalues $\lambda_1,\lambda_2,\lambda_3\in\rel$, and let $\{x_1,x_2,y_2,x_3,y_3\}$ be an orthonormal basis of $(\m,\met)$ such that:
	\begin{equation*}
		\m^\h=\m_1=\ker\big(\rico-\lambda_1\id_\m\big)=\rel x_1,\et [\h,\m]=\m_2\oplus\m_3,
	\end{equation*}
	where $\m_i=\ker\big(\rico-\lambda_i\id_\m\big)=\Span_\rel\{x_i,y_i\}$ for $i\in\{2,3\}$. Since $(G/H,\met)$ is Ricci-Codazzi, using $\big[\m^\h,[\h,\m]\big]_\m\subseteq [\h,\m]$, Corollary \ref{eigenspacesstable}, and \eqref{CodazIeqlJ}, we can write:
	\begin{equation}\label{bracketm}
		\left\{
		\begin{aligned}
			[x_1,x_2]_\m\,\,&=&\varepsilon_1y_2+\mu_1x_3+\gamma_1y_3;\quad\\
			[x_1,y_2]_\m\,\,&=&-\varepsilon_1x_2+\mu_2x_3+\gamma_2y_3;\quad\\
			[x_1,x_3]_\m\,\,&=&\mu_3x_2+\gamma_3y_2+\varepsilon_2y_3;\quad\\
			[x_1,y_3]_\m\,\,&=&\mu_4x_2+\gamma_4y_2-\varepsilon_2x_3;\quad\\
			[x_2,x_3]_\m\,\,&=&\mu_5x_1+\varepsilon_3y_2+\varepsilon_4y_3;\quad\\
			[x_2,y_3]_\m\,\,&=&\gamma_5x_1+\varepsilon_5y_2-\varepsilon_4x_3;\quad\\
			[y_2,x_3]_\m\,\,&=&\gamma_6x_1-\varepsilon_3x_2+\varepsilon_6y_3;\quad\\
			[y_2,y_3]_\m\,\,&=&\mu_6x_1-\varepsilon_5x_2-\varepsilon_6x_3,\quad
		\end{aligned}
		\right.
	\end{equation}
	where $\varepsilon_i,\mu_i,\gamma_i\in\rel$ for $i\in\{1,\ldots,6\}$. From  Corollary \ref{CodazziParallel} one can conclude that $\n\ric=0$ if and only if $\mu_i=\gamma_i=0$ for all $i\in\{1,\ldots,6\}$. Combining \eqref{EquivCodazzi} with \eqref{bracketm}, we obtain:
	\begin{equation}\label{ReductionCoef}
		\left\{
		\begin{aligned}
			\mu_3\,&=&-\tfrac{\Delta_{13}^2}{\Delta_{12}^2}\mu_1;\,\\
			\mu_4\,&=&-\tfrac{\Delta_{13}^2}{\Delta_{12}^2}\gamma_1;\,\\
			\mu_5\,&=&\tfrac{\Delta_{23}^2}{\Delta_{12}^2}\mu_1;\,\\
			\mu_6\,&=&\tfrac{\Delta_{23}^2}{\Delta_{12}^2}\gamma_2;\,
		\end{aligned}
		\right.\qquad\et\qquad
		\left\{
		\begin{aligned}
			\gamma_3\,&=&-\tfrac{\Delta_{13}^2}{\Delta_{12}^2}\mu_2;\,\\
			\gamma_4\,&=&-\tfrac{\Delta_{13}^2}{\Delta_{12}^2}\gamma_2;\,\\
			\gamma_5\,&=&\tfrac{\Delta_{23}^2}{\Delta_{12}^2}\gamma_1;\,\\
			\gamma_6\,&=&\tfrac{\Delta_{23}^2}{\Delta_{12}^2}\mu_2,\,
		\end{aligned}
		\right.
	\end{equation}
	where $\Delta_{ij}:=\lambda_i-\lambda_j$ for $i,j\in\{1,2,3\}$. Thus, we must show that $\mu_i=\gamma_i=0$ for $i\in\{1,2\}$. Let us compute the different products $\Ll_{x_i}(x_j)$, $\Ll_{x_i}(y_j)$, $\Ll_{y_i}(x_j)$, and $\Ll_{y_i}(y_j)$. Note that $[\m_2,\m_2]_\m=[\m_3,\m_3]_\m=\{0\}$ is a direct consequence of Corollary \ref{eigenspacesstable} and the second point of Proposition \ref{propertiesreductive}, moreover, using \eqref{LeviCivita}, we deduce that $\Ll(\m_2,\m_2)=\Ll(\m_3,\m_3)=\{0\}$, meaning
	\begin{equation*}
		\Ll_{x_i}(x_i)=\Ll_{y_j}(y_j)=\Ll_{x_i}(y_i)=\Ll_{y_j}(x_j)=0,
	\end{equation*}
	for all $i,j\in\{2,3\}$. We also have $\Ll_{x_1}(x_1)=0$. Similarly, a straightforward computation, using \eqref{LeviCivita} and \eqref{bracketm}, gives
	\begin{equation}\label{Civitam}
		\left\{
		\begin{aligned}
			\Ll_{x_1}(x_2)\,&=&\varepsilon_1y_2-\tfrac{1}{2}(-\mu_1+\mu_3+\mu_5)x_3-\tfrac{1}{2}(-\gamma_1+\gamma_5+\mu_4)y_3;\hfill\\
			\Ll_{x_1}(x_3)\,&=&\tfrac{1}{2}(\mu_3-\mu_1+\mu_5)x_2+\tfrac{1}{2}(\gamma_3-\mu_2+\gamma_6)y_2+\varepsilon_2y_3;\quad\,\,\\
			\Ll_{x_1}(y_2)\,&=&-\varepsilon_1x_2-\tfrac{1}{2}(-\mu_2+\gamma_6+\gamma_3)x_3-\tfrac{1}{2}(-\gamma_2+\mu_6+\gamma_4)y_3;\hfill\\
			\Ll_{x_1}(y_3)\,&=&\tfrac{1}{2}(\mu_4+\gamma_5-\gamma_1)x_2+\tfrac{1}{2}(\gamma_4+\mu_6-\gamma_2)y_2-\varepsilon_2x_3;\quad\,\,\\
			\Ll_{x_2}(x_1)\,&=&-\tfrac{1}{2}(\mu_1+\mu_3+\mu_5)x_3-\tfrac{1}{2}(\gamma_1+\gamma_5+\mu_4)y_3;\qquad\quad\quad\,\,\\
			\Ll_{x_2}(x_3)\,&=&\tfrac{1}{2}(\mu_1+\mu_3+\mu_5)x_1+\varepsilon_4y_3;\qquad\qquad\qquad\qquad\qquad\\
			\Ll_{x_2}(y_3)\,&=&\tfrac{1}{2}(\mu_4+\gamma_1+\gamma_5)x_1-\varepsilon_4x_3.\qquad\qquad\qquad\qquad\qquad
		\end{aligned}
		\right.
	\end{equation}
	Moreover,
	\begin{equation}
		\left\{
		\begin{aligned}
			\Ll_{y_2}(x_1)\,&=&-\tfrac{1}{2}(\mu_2+\gamma_6+\gamma_3)x_3-\tfrac{1}{2}(\gamma_2+\mu_6+\gamma_4)y_3;\qquad\quad\quad\,\,\\
			\Ll_{y_2}(x_3)\,&=&\tfrac{1}{2}(\mu_2+\gamma_6+\gamma_3)x_1+\varepsilon_6y_3;\qquad\qquad\qquad\qquad\qquad\,\,\\
			\Ll_{y_2}(y_3)\,&=&\tfrac{1}{2}(\gamma_2+\gamma_4+\mu_6)x_1-\varepsilon_6x_3;\qquad\qquad\qquad\qquad\qquad\,\\
			\Ll_{x_3}(x_1)\,&=&-\tfrac{1}{2}(\mu_3+\mu_1-\mu_5)x_2-\tfrac{1}{2}(\gamma_3+\mu_2-\gamma_6)y_2;\qquad\quad\quad\,\,\\
			\Ll_{x_3}(x_2)\,&=&\tfrac{1}{2}(\mu_1+\mu_3-\mu_5)x_1-\varepsilon_3y_2;\qquad\qquad\qquad\qquad\qquad\\
			\Ll_{x_3}(y_2)\,&=&\tfrac{1}{2}(\mu_2+\gamma_3-\gamma_6)x_1+\varepsilon_3x_2;\qquad\qquad\qquad\qquad\qquad\\
			\Ll_{y_3}(x_1)\,&=&-\tfrac{1}{2}(\mu_4+\gamma_1-\gamma_5)x_2-\tfrac{1}{2}(\gamma_4+\gamma_2-\mu_6)y_2;\qquad\quad\quad\,\,\\
			\Ll_{y_3}(x_2)\,&=&\tfrac{1}{2}(\mu_4+\gamma_1-\gamma_5)x_1-\varepsilon_5y_2,\qquad\qquad\qquad\qquad\qquad\\
			\Ll_{y_3}(y_2)\,&=&\tfrac{1}{2}(\gamma_2+\gamma_4-\mu_6)x_1+\varepsilon_5x_2.\qquad\qquad\qquad\qquad\qquad\\
		\end{aligned}
		\right.
	\end{equation}
	On the other hand, evaluating $\ric(x_k,x_k)$ and $\ric(y_k,y_k)$ yields explicit expressions for the eigenvalues $\lambda_1,\lambda_2,\lambda_3$ in terms of the parameters $\mu_i$ and $\gamma_j$, i.e.,
	\begin{equation}\label{Lambdas}
		\left\{
		\begin{aligned}
			\lambda_1\,=\,\,\ric(x_1,x_1)\,&=&\tfrac{1}{2}\Big\{(\mu_1+\mu_3)^2+(\mu_2+\gamma_3)^2+(\mu_4+\gamma_1)^2+(\gamma_2+\gamma_4)^2-\mu_5^2-\mu_6^2-\gamma_5^2-\gamma_6^2\Big\};\\
			\lambda_2\,=\,\,\ric(x_2,x_2)\,&=&\tfrac{1}{2}\Big\{(\mu_1-\mu_5)^2+(\gamma_1-\gamma_5)^2-\mu_3^2-\mu_4^2\Big\}-\alpha([x_2,y_2]);\qquad\qquad\qquad\qquad\qquad\,\\
			\lambda_2\,=\,\,\ric(y_2,y_2)\,&=&\tfrac{1}{2}\Big\{(\mu_2-\gamma_6)^2+(\gamma_2-\mu_6)^2-\gamma_3^2-\gamma_4^2\Big\}-\alpha([x_2,y_2]);\qquad\qquad\qquad\qquad\qquad\,\\
			\lambda_3\,=\,\,\ric(x_3,x_3)\,&=&\tfrac{1}{2}\Big\{(\mu_3+\mu_5)^2+(\gamma_3+\gamma_6)^2-\mu_1^2-\mu_2^2\Big\}-\beta([x_3,y_3]);\qquad\qquad\qquad\qquad\qquad\,\\
			\lambda_3\,=\,\,\ric(y_3,y_3)\,&=&\tfrac{1}{2}\Big\{(\mu_4+\gamma_5)^2+(\mu_6+\gamma_4)^2-\gamma_1^2-\gamma_2^2\Big\}-\beta([x_3,y_3]),\qquad\qquad\qquad\qquad\qquad\,\\
		\end{aligned}
		\right.
	\end{equation}
	where $\alpha,\beta\in\h^*\backslash\{0\}$ such that:
	\begin{equation*}
		[u,x_2]=\alpha(u)y_2,\quad [u,y_2]=-\alpha(u)x_2,\quad [u,x_3]=\alpha(u)y_3,\et [u,y_3]=-\alpha(u)x_3,
	\end{equation*}
	for all $u\in\h$ (see Corollary \ref{CorDim2}). In addition, we have:
	\begin{equation}\label{Ric=0}
		\left\{
		\begin{aligned}
			0\,=\,\,\ric(x_2,y_2)\,&=&\tfrac{1}{2}\Big\{(\gamma_2-\mu_6)(\gamma_1-\gamma_5)+(\mu_1-\mu_5)(\mu_2-\gamma_6)-\mu_3\gamma_3-\mu_4\gamma_4\Big\};\qquad\qquad\qquad\qquad\qquad\,\\
			0\,=\,\,\ric(x_3,y_3)\,&=&\tfrac{1}{2}\Big\{(\mu_3+\mu_5)(\mu_4+\gamma_5)+(\mu_6+\gamma_4)(\gamma_3+\gamma_6)-\mu_1\gamma_1-\mu_2\gamma_2\Big\}.\qquad\qquad\qquad\qquad\qquad\,
		\end{aligned}
		\right.
	\end{equation}
	Substituting \eqref{ReductionCoef} in \eqref{Ric=0}, we get:
	\begin{equation}\label{systric=0}
		\left\{
		\begin{aligned}
			\gamma_1\gamma_2+\mu_1\mu_2\,&=&0;\\
			\mu_1\gamma_1+\mu_2\gamma_2\,&=&0.
		\end{aligned}
		\right.
	\end{equation}
	We claim that $\gamma_1^2=\mu_2^2$ and $\gamma_2^2=\mu_1^2$. In fact, assume by contradiction that $\gamma_1^2\neq\mu_2^2$, then from \eqref{systric=0} we deduce that $\gamma_2=\mu_1=0$. In particular, by \eqref{ReductionCoef} we have $\mu_3=\mu_5=\mu_6=\gamma_4=0$. Therefore, using \eqref{Lambdas} and substituting in the equation $\ric(y_3,y_3)-\ric(x_3,x_3)=0$, we get that $\gamma_1^2=\mu_2^2$, which is impossible. Hence, $\gamma_1^2=\mu_2^2$. An analogous argument, using $\ric(x_2,x_2)-\ric(y_2,y_2)=0$, gives that $\gamma_2^2=\mu_1^2$. Thus, to prove that $(G/H,\met)$ is Ricci-parallel it suffices to show that $(\mu_1,\mu_2)=(0,0)$ or equivalently $(\gamma_1,\gamma_2)=(0,0)$. If $(\varepsilon_3,\varepsilon_4,\varepsilon_5,\varepsilon_6)\in\rel^4\backslash\{(0,0,0,0)\}$, then a direct computation yields
	\begin{equation*}
		\left\{
		\begin{aligned}
			\langle\R(x_1,x_3)y_2,x_3\rangle\,\,&=&\frac{\mu_1\varepsilon_3(\Delta_{12}^2+\Delta_{13}^2-\Delta_{23}^2)}{2\Delta_{12}^2}\,\,=\,\,\frac{\mu_1\varepsilon_3\Delta_{13}}{\Delta_{12}};\quad\\
			\langle\R(x_1,x_3)x_2,x_3\rangle\,\,&=&-\frac{\mu_2\varepsilon_3(\Delta_{12}^2+\Delta_{13}^2-\Delta_{23}^2)}{2\Delta_{12}^2}\,\,=\,\,-\frac{\mu_2\varepsilon_3\Delta_{13}}{\Delta_{12}};\\
			\langle\R(x_1,x_2)y_3,x_2\rangle\,\,&=&\frac{\mu_1\varepsilon_4(\Delta_{12}^2+\Delta_{13}^2-\Delta_{23}^2)}{2\Delta_{12}^2}\,\,=\,\,\frac{\mu_1\varepsilon_4\Delta_{13}}{\Delta_{12}};\quad\\
			\langle\R(x_1,x_2)y_3,y_2\rangle\,\,&=&\frac{\mu_2\varepsilon_4(\Delta_{12}^2+\Delta_{13}^2-\Delta_{23}^2)}{2\Delta_{12}^2}\,\,=\,\,\frac{\mu_2\varepsilon_4\Delta_{13}}{\Delta_{12}};\quad\\
			\langle\R(y_2,y_3)x_1,x_3\rangle\,\,&=&-\frac{\mu_1\varepsilon_5(\Delta_{12}^2+\Delta_{23}^2-\Delta_{13}^2)}{2\Delta_{12}^2}\,\,=\,\,-\frac{\mu_1\varepsilon_5\Delta_{23}}{\Delta_{12}};\\
			\langle\R(x_2,y_3)x_1,x_3\rangle\,\,&=&\frac{\mu_2\varepsilon_5(\Delta_{12}^2+\Delta_{23}^2-\Delta_{13}^2)}{2\Delta_{12}^2}\,\,=\,\,\frac{\mu_2\varepsilon_5\Delta_{23}}{\Delta_{12}};\quad\\
			\langle\R(y_2,y_3)x_1,x_2\rangle\,\,&=&\frac{\mu_1\varepsilon_6(\Delta_{13}^2+\Delta_{23}^2-\Delta_{12}^2)}{2\Delta_{12}^2}\,\,=\,\,\frac{\mu_1\varepsilon_6\Delta_{13}\Delta_{23}}{\Delta^2_{12}};\\
			\langle\R(y_2,y_3)x_1,y_2\rangle\,\,&=&\frac{\mu_2\varepsilon_6(\Delta_{13}^2+\Delta_{23}^2-\Delta_{12}^2)}{2\Delta_{12}^2}\,\,=\,\,\frac{\mu_2\varepsilon_6\Delta_{13}\Delta_{23}}{\Delta^2_{12}}.
		\end{aligned}
		\right.
	\end{equation*}
	Thus, using \eqref{Curv=0}, we obtain $\mu_1=\mu_2=0$ and the conclusion follows. In what follows, assume that $\varepsilon_3=\varepsilon_4=\varepsilon_5=\varepsilon_6=0$. Then, \eqref{bracketm} becomes:
	\begin{equation}\label{bracketmp}
		\left\{
		\begin{aligned}
			[x_1,x_2]_\m\,\,&=&\varepsilon_1y_2+\mu_1x_3+\gamma_1y_3;\quad\\
			[x_1,y_2]_\m\,\,&=&-\varepsilon_1x_2+\mu_2x_3+\gamma_2y_3;\quad\\
			[x_1,x_3]_\m\,\,&=&\mu_3x_2+\gamma_3y_2+\varepsilon_2y_3;\quad\\
			[x_1,y_3]_\m\,\,&=&\mu_4x_2+\gamma_4y_2-\varepsilon_2x_3;\quad\\
			[x_2,x_3]_\m\,\,&=&\mu_5x_1;\qquad\qquad\qquad\quad\\
			[x_2,y_3]_\m\,\,&=&\gamma_5x_1;\qquad\qquad\qquad\quad\\
			[y_2,x_3]_\m\,\,&=&\gamma_6x_1;\qquad\qquad\qquad\quad\\
			[y_2,y_3]_\m\,\,&=&\mu_6x_1.\qquad\qquad\qquad\quad
		\end{aligned}
		\right.
	\end{equation}
	In particular, $[\m_2,\m_3]_\m\subseteq\m_1$. Now, assume by contradiction that $(\gamma_1,\gamma_2)\neq (0,0)$ and therefore $(\mu_1,\mu_2)\neq (0,0)$. For any $u\in\h$, we compute:
	\begin{eqnarray*}
		0&=&\langle\big[u,[x_2,x_3]_\m\big],x_1\rangle\\
		&=&\langle\big[[u,x_2],x_3\big]_\m,x_1\rangle+\langle\big[x_2,[u,x_3]\big]_\m,x_1\rangle\\
		&\stackrel{\rm\eqref{bracketmp}}{=}&\frac{\Delta_{23}^2}{\Delta_{12}^2}\Big(\mu_2\alpha(u)+\gamma_1\beta(u)\Big).
	\end{eqnarray*}
	Analogously, we obtain:
	\begin{equation*}
		0\,=\,\langle\big[u,[y_2,x_3]_\m\big],x_1\rangle\,=\,\frac{\Delta_{23}^2}{\Delta_{12}^2}\Big(\gamma_2\beta(u)-\mu_1\alpha(u)\Big).
	\end{equation*}
	As a result,
	\begin{equation}\label{AlphaBeta}
		\left\{
		\begin{aligned}
			\mu_2\alpha(u)\,&=&-\gamma_1\beta(u);\qquad\qquad\quad\,\\
			\gamma_2\beta(u)\,&=&\mu_1\alpha(u),\qquad\forall\,u\in\h.
		\end{aligned}
		\right.
	\end{equation}
	Let us now give explicit expressions for $\alpha([x_2,y_2])$ and $\beta([x_3,y_3])$ in terms of the eigenvalues $\lambda_1,\lambda_2,\lambda_3$. First, we compute $\alpha([x_3,y_3])$ and $\beta([x_2,y_2])$:
	\begin{equation*}0=\langle\big[[x_2,x_3],y_3\big]_\m+\big[[y_3,x_2],x_3\big]_\m+\big[[x_3,y_3],x_2\big]_\m,y_2\rangle=\big(\alpha([x_3,y_3])+\mu_5\gamma_4-\gamma_3\gamma_5\big).
	\end{equation*}
	Thus, 
	\begin{equation}\label{alphx3y3}
		\alpha([x_3,y_3])=-\frac{\Delta^2_{13}\Delta^2_{23}}{\Delta^4_{12}}\big(\mu_2\gamma_1-\mu_1\gamma_2\big).
	\end{equation}
	Likewise,
	\begin{equation*}
		0=\langle\big[[x_3,x_2],y_2\big]_\m+\big[[y_2,x_3],x_2\big]_\m+\big[[x_2,y_2],x_3\big]_\m,y_3\rangle=\big(\beta([x_2,y_2])-\mu_5\gamma_2+\gamma_1\gamma_6\big).
	\end{equation*}
	Hence,
	\begin{equation}\label{betx2y2}
		\beta([x_2,y_2])=\frac{\Delta^2_{23}}{\Delta^2_{12}}\big(\mu_1\gamma_2-\mu_2\gamma_1\big).
	\end{equation}
	If $\gamma_1\neq 0$ (equivalently $\mu_2\neq 0$), then from \eqref{systric=0}, the first equality in \eqref{AlphaBeta}, \eqref{alphx3y3}, and \eqref{betx2y2}, we deduce that:
	\begin{equation*}
		\left\{
		\begin{aligned}
			\alpha([x_2,y_2])\,&=&-\frac{\gamma_1}{\mu_2}\beta([x_2,x_2])\,=\,\frac{\Delta^2_{23}}{\Delta^2_{12}}(\gamma_1^2+\gamma_2^2);\qquad\\
			\beta([x_3,y_3])\,&=&-\frac{\mu_2}{\gamma_1}\alpha([x_3,x_3])\,=\,\frac{\Delta^2_{13}\Delta^2_{23}}{\Delta^4_{12}}(\gamma_1^2+\gamma_2^2).\,\,
		\end{aligned}
		\right.
	\end{equation*}
	If $\gamma_2\neq 0$ (equivalently $\mu_1\neq 0$), then from \eqref{systric=0}, the second equality in \eqref{AlphaBeta}, \eqref{alphx3y3}, and \eqref{betx2y2}, we get that:
	\begin{equation}\label{alx2y2bex3y3}
		\left\{
		\begin{aligned}
			\alpha([x_2,y_2])\,&=&\frac{\Delta^2_{23}}{\Delta^2_{12}}(\gamma_1^2+\gamma_2^2);\qquad\\
			\beta([x_3,y_3])\,&=&\frac{\Delta^2_{13}\Delta^2_{23}}{\Delta^4_{12}}(\gamma_1^2+\gamma_2^2).\,\,
		\end{aligned}
		\right.
	\end{equation}
	Regardless of either $\gamma_1\neq 0$ or $\gamma_2\neq 0$, we have shown that \eqref{alx2y2bex3y3} always holds. Next, using that $\mu_1^2=\gamma_2^2$, $\mu_2^2=\gamma_1^2$, and \eqref{Lambdas}, we have:
	\begin{eqnarray}\label{lambda1g1g2}
		\lambda_1&=&\tfrac{1}{2}\Big\{(\mu_1+\mu_3)^2+(\mu_2+\gamma_3)^2+(\mu_4+\gamma_1)^2+(\gamma_2+\gamma_4)^2-\mu_5^2-\mu_6^2-\gamma_5^2-\gamma_6^2\Big\}\nonumber\\
		&=&4\frac{\Delta_{13}\Delta^2_{23}}{\Delta^3_{12}}\big(\gamma^2_1+\gamma_2^2\big).
	\end{eqnarray}
	In particular, $\lambda_1\neq 0$. By a similar computation, using \eqref{alx2y2bex3y3}, we obtain:
	\begin{equation}\label{lambda2lambda3g1g2}
		\lambda_2=-\frac{\Delta_{23}}{\Delta_{12}^3}\Big(2\Delta^2_{13}+\Delta_{12}\Delta_{23}\Big)\big(\gamma^2_1+\gamma_2^2\big),\et \lambda_3=\frac{\Delta_{13}\Delta_{23}}{\Delta_{12}^4}\Big(2\Delta^2_{12}-\Delta_{13}\Delta_{23}\Big)\big(\gamma^2_1+\gamma_2^2\big).
	\end{equation}
	From equality \eqref{lambda1g1g2}, we get that:
	\begin{equation}
		\gamma^2_1+\gamma_2^2=\frac{\lambda_1\Delta_{12}^3}{4\Delta_{13}\Delta^2_{23}}
	\end{equation}
	Hence, the equalities in \eqref{lambda2lambda3g1g2} show that:
	\begin{equation*}
		\lambda_2=-\frac{\lambda_1}{4}\Big(\frac{2\Delta_{13}}{\Delta_{23}}+\frac{\Delta_{12}}{\Delta_{13}}\Big),\et \lambda_3=\frac{\lambda_1}{4}\Big(\frac{2\Delta_{12}}{\Delta_{23}}-\frac{\Delta_{13}}{\Delta_{12}}\Big).
	\end{equation*}
	Consequently,
	\begin{equation*}
		\lambda_1-\lambda_2=\frac{\lambda_1}{4}\Big(4+\frac{2\Delta_{13}}{\Delta_{23}}+\frac{\Delta_{12}}{\Delta_{13}}\Big)=\frac{\lambda_1\Big(4\Delta_{13}\Delta_{23}+2\Delta_{13}^2+\Delta_{12}\Delta_{23}\Big)}{4\Delta_{13}\Delta_{23}}.
	\end{equation*}
	Hence,
	\begin{equation}\label{delt12}
		4\Delta_{12}\Delta_{13}\Delta_{23}=\lambda_1\Big(4\Delta_{13}\Delta_{23}+2\Delta_{13}^2+\Delta_{12}\Delta_{23}\Big).
	\end{equation}
	On the other hand, 
	\begin{equation*}
		\lambda_1-\lambda_3=\frac{\lambda_1}{4}\Big(4-\frac{2\Delta_{12}}{\Delta_{23}}+\frac{\Delta_{13}}{\Delta_{12}}\Big)=\frac{\lambda_1\Big(4\Delta_{12}\Delta_{23}-2\Delta_{12}^2+\Delta_{13}\Delta_{23}\Big)}{4\Delta_{12}\Delta_{23}}.
	\end{equation*}
	Therefore,
	\begin{equation}\label{delt13}
		4\Delta_{12}\Delta_{13}\Delta_{23}=\lambda_1\Big(4\Delta_{12}\Delta_{23}-2\Delta_{12}^2+\Delta_{13}\Delta_{23}\Big).
	\end{equation}
	Now, using the fact that $\lambda_1\neq 0$ with $\eqref{delt12}-\eqref{delt13}=0$, we get:
	\begin{equation*}
		2\Delta_{12}^2+2\Delta_{13}^2+3\Delta_{23}^2=0.
	\end{equation*}
	However, this contradict the fact that $\lambda_1,\lambda_2,$ and $\lambda_3$ are mutually distinct. As a result, we obtain that $\gamma_1=\gamma_2=0$, which shows that $(G/H,\met)$ is Ricci-parallel.
	
	\subsection{The case where $\m^\h$ is $2$-dimensional}
	
	In this case $[\h,\m]$ is $3$-dimensional, and hence it is an irreducible $\h$-submodule of $\m$ by Proposition \ref{propertiesreductive}. It follows from \eqref{RicStable} that $\rico_{[\h,\m]}=\lambda\id_{[\h,\m]}$ with $\lambda\in\rel$. As a result, $\rico:\m\to\m$ has at most three eigenvalues. If $\rico:\m\to\m$ has two eigenvalues, then the Ricci parallelness follows from Proposition \ref{Twoeigenval}. If $\rico:\m\to\m$ has three eigenvalues, then using the two inclusions $[\m^\h,\m^\h]_\m\subseteq\m^\h$ and $\big[\m^\h,[\h,\m]\big]_\m\subseteq [\h,\m]$, the result follows from Corollary \ref{CodazziParallel}.

	\subsection{The case where $\m^\h$ is $3$-dimensional}
	
	In this case $[\h,\m]$ is $2$-dimensional and therefore by virtue of Proposition \ref{propertiesreductive} and \eqref{RicStable}, we have $\rico_{[\h,\m]}=\lambda\id_{[\h,\m]}$ with $\lambda\in\rel$. Thus, the number of eigenvalues of  $\rico:\m\to\m$ in this case is at most four. If $\rico:\m\to\m$ has less than three eigenvalues then the claim follows from Propositions \ref{EinstParallel} and \ref{Twoeigenval}.  On the other hand, if $\rico$ has exactly three eigenvalues such that $\m^\h$ splits as the direct sum of two eigenspaces, then Lemmas \ref{LemInvVect}, \ref{lemDecompM} together with Corollary \ref{CodazziParallel}, implies that $(G/H,\met)$ is Ricci-parallel. Therefore the only remaining cases to consider are the following:
	\begin{itemize}
		\item[$\bullet$] $\rico:\m\to\m$ has four distinct eigenvalues $\lambda_1,\lambda_2,\lambda_3$ and $\lambda_4$.
		\item[$\bullet$] $\rico:\m\to\m$ has three distinct eigenvalues $\lambda_1,\lambda_2,\lambda_3$ such that $[\h,\m]\subsetneqq\m_3$ and $\m_3\cap\m^\h\neq \{0\}$, where $\m_3:=\ker\big(\rico-\lambda_3\id_\m\big)$.
	\end{itemize}
	First, assume that $\rico:\m\to\m$ has four mutually distinct eigenvalues $\lambda_1,\lambda_2,\lambda_3$ and $\lambda_4$. Let $\{x_1,x_2,x_3,x_4,y_4\}$ be an orthonormal basis of $(\m,\met)$ such that:
	\begin{equation*}
		\m^\h=\Span_\rel\{x_1,x_2,x_3\},\quad [\h,\m]=\Span_\rel\{x_4,y_4\},\quad \rico(x_i)=\lambda_i x_i,\et \rico(y_4)=\lambda_4y_4,
	\end{equation*} 
	for $i\in\{1,2,3,4\}$. Since $[\h,\m]$ is an eigenspace of $\rico:\m\to\m$, then combining Corollary \ref{CorDim2} with Corollary \ref{eigenspacesstable} one can deduce that $[x_4,y_4]_\m=0$. Moreover, using \eqref{CodazIeqlJ} and the fact that $\m^\h$ is a subalgebra of $(\m,\cro_\m)$, we can find three real numbers $\mu_1,\mu_2,\mu_3\in\rel$ such that:
	\begin{equation*}
		[x_1,x_2]_\m=\mu_3x_3,\quad [x_1,x_3]_\m=\mu_2x_2,\et [x_2,x_3]_\m=\mu_1x_1.
	\end{equation*}
	From \eqref{EquivCodazzi}, the coefficients $\mu_1,\mu_2$, and $\mu_3$ satisfy
	\begin{eqnarray}
		\mu_2(\lambda_1-\lambda_2)^2+\mu_3(\lambda_1-\lambda_3)^2&=&0;\label{equation1dim3}\\
		\mu_1(\lambda_2-\lambda_1)^2-\mu_3(\lambda_2-\lambda_3)^2&=&0.\label{equation2dim3}
	\end{eqnarray}
	Further, using $(\lambda_2-\lambda_3)^2\eqref{equation1dim3}+(\lambda_1-\lambda_3)^2\eqref{equation2dim3}=0$, we obtain:
	$$\mu_1(\lambda_1-\lambda_3)^2+\mu_2(\lambda_2-\lambda_3)^2=0.$$
	In particular, by Corollary \ref{CodazziParallel}, we have $\n\ric=0$ if and only if $\mu_1\mu_2\mu_3=0$.
	Assume by contradiction that $\mu_3\neq 0$ and define $\Delta_{ij}:=\lambda_i-\lambda_j$ for $i,j\in\{1,2,3\}$. Then from \eqref{equation1dim3} and \eqref{equation2dim3}, we have:
	\begin{equation}\label{mu1mu2dim3}
		\mu_1=\frac{\Delta_{23}^2}{\Delta_{12}^2}\mu_3,\et \mu_2=-\frac{\Delta_{13}^2}{\Delta_{12}^2}\mu_3.
	\end{equation}
	On the other hand, a straightforward computation yields:
	\begin{equation}\label{la1dim3}
		\lambda_1=\ric(x_1,x_1)=\frac{(\mu_2+\mu_3)^2-\mu_1^2}{2}.
	\end{equation}
	Similarly, we get:
	\begin{equation}\label{la2la3dim3}
		\lambda_2=\frac{(\mu_1-\mu_3)^2-\mu_2^2}{2},\et \lambda_3=\frac{(\mu_1+\mu_2)^2-\mu_3^2}{2}.
	\end{equation}
	Substituting \eqref{mu1mu2dim3} in \eqref{la1dim3} and \eqref{la2la3dim3}, we obtain:
	\begin{equation*}
		\lambda_1=\frac{2\Delta_{13}\Delta_{23}^2}{\Delta_{12}^3}\mu_3^2,\quad\lambda_2=-\frac{2\Delta_{13}^2\Delta_{23}}{\Delta_{12}^3}\mu_3^2,\et \lambda_3=\frac{2\Delta_{13}\Delta_{23}}{\Delta_{12}^2}\mu_3^2.
	\end{equation*}
	Substituting again by $\mu_3^2=\tfrac{\Delta_{12}^2}{2\Delta_{13}\Delta_{23}}\lambda_3$ in $\lambda_1$ and $\lambda_2$, we have:
	\begin{equation*}
		\lambda_1=\frac{\Delta_{23}}{\Delta_{12}}\lambda_3,\et \lambda_2=-\frac{\Delta_{13}}{\Delta_{12}}\lambda_3.
	\end{equation*}
	Thus, 
	\begin{equation*}
		\Delta_{13}=\frac{\Delta_{23}-\Delta_{12}}{\Delta_{12}}\lambda_3,\et \Delta_{23}=-\frac{\Delta_{13}+\Delta_{12}}{\Delta_{12}}\lambda_3.
	\end{equation*}
	It follows that,
	\begin{equation*}
		\Delta_{12}\Delta_{13}\Delta_{23}=\big(\Delta_{23}-\Delta_{12}\big)\Delta_{23}\lambda_3,\et \Delta_{12}\Delta_{23}\Delta_{13}=-\big(\Delta_{13}+\Delta_{12}\big)\Delta_{13}\lambda_3.
	\end{equation*}
	As a result, since $\lambda_3\neq 0$, we get:
	\begin{equation*}
		\Delta_{12}^2+\Delta_{13}^2+\Delta_{23}^2=0,
	\end{equation*}
	which is a contradiction. Hence, $\mu_1=\mu_2=\mu_3=0$, and therefore $(G/H,\met)$ is Ricci-parallel.
	
	Now, suppose that $\rico:\m\to\m$ has three mutually distinct eigenvalues $\lambda_1,\lambda_2,\lambda_3$ such that:
	\begin{equation*}
		\m_1=\rel x_1,\,\,\m_2=\rel x_2,\,\,\m_3=\rel x_3\oplus [\h,\m],\,\,\m^\h=\Span_\rel\{x_1,x_2,x_3\},\,\,\textnormal{and }\,\, [\h,\m]=\Span_\rel\{y_3,z_3\},
	\end{equation*}
	where $\m_i:=\ker\big(\rico-\lambda_i\id_\m\big)$ for $i\in\{1,2,3\}$, and $\{x_1,x_2,x_3,y_3,z_3\}$ is an orthonormal basis of $(\m,\met)$. By an argument similar to that of the previous case, we can write:
	\begin{equation*}
		[x_1,x_2]_\m=\mu_3x_3,\quad [x_1,x_3]_\m=\mu_2x_2,\et [x_2,x_3]_\m=\mu_1x_1,
	\end{equation*}
	for $\mu_1,\mu_2,\mu_3\in\rel$. It is easy to verify that \eqref{equation1dim3} and \eqref{equation2dim3} both hold in this case. Therefore, to prove the claim, it suffices to show that $\mu_3=0$. The fact that $\LS(\m_3,\m_3)\subseteq \m_3$ implies that $\LS\big(\m_3,\m_1\oplus\m_2\big)\subseteq\m_1\oplus\m_2$. Moreover, since $[\h,\m]\subseteq\m_3$, it follows that $\LS\big([\h,\m],\m_1\oplus\m_2\big)\subseteq\m_1\oplus\m_2$. However, we also have $\LS\big([\h,\m],\m^\h\big)\subseteq[\h,\m]$, and therefore $\LS\big([\h,\m],\m_1\oplus\m_2\big)\subseteq[\h,\m]$. Combining both facts shows that $\LS\big([\h,\m],\m_1\oplus\m_2\big)=\{0\}$. A straightforward computation of $\ric(x_1,x_1)$ and $\ric(x_2,x_2)$, using the previous discussion, shows:
	\begin{equation}\label{Lamd1andLamd2}
		\lambda_1=\frac{(\mu_2+\mu_3)^2-\mu_1^2}{2},\et \lambda_2=\frac{(\mu_1-\mu_3)^2-\mu_2^2}{2}.
	\end{equation}
	On the other hand, since $\big[[\h,\m],[\h,\m]\big]_\m$ is a $1$-dimensional $\h$-submodule of $\m$, it is contained in $\m^\h$ (see Proposition \ref{propertiesreductive}). Therefore, using Corollary \ref{eigenspacesstable}, we can write $[y_3,z_3]_\m=\nu x_3$ with $\nu\in\rel$. Further, since $[\m^\h,[\h,\m]]_\m\subseteq [\h,\m]$, we have:
	\begin{equation*}
		[x_3,y_3]_\m=\gamma y_3+\eta z_3,\et [x_3,z_3]_\m=\varepsilon y_3+\xi z_3,
	\end{equation*}
	where $\gamma,\eta,\varepsilon,\xi\in\rel$. In addition, for any $u\in\h$, we can write $[u,y_3]=\alpha(u)z_3$ and $[u,z_3]=-\alpha(u)y_3$, where $\alpha\in\h^*\backslash\{0\}$. Furthermore, since $\LS(\m^\h,[\h,\m])\subseteq [\h,\m]$, $\LS\big([\h,\m],\m_1\oplus\m_2\big)=\{0\}$, and $\LS_x\in\mathfrak{so}(\m,\met)$ for any $x\in\m$, we get that:
	\begin{equation*}
		\left\{
		\begin{aligned}
			[x_1,y_3]_\m\,&=&x_1\cdot y_3\,\,=\,\,az_3;\qquad\quad\\
			[x_1,z_3]_\m\,&=&x_1\cdot z_3\,\,=\,\,-ay_3;\qquad\,\\
			[x_2,y_3]_\m\,&=&x_2\cdot y_3\,\,=\,\,bz_3;\qquad\quad\\
			[x_2,z_3]_\m\,&=&x_2\cdot z_3\,\,=\,\,-by_3.\qquad\,\,
		\end{aligned}
		\right.
	\end{equation*}
	We compute:
	\begin{equation}\label{Jac1}
		\big[[x_1,x_2]_\m,y_3\big]_\m=\big[[x_1,y_3]_\m,x_2\big]_\m+\big[x_1,[x_2,y_3]_\m\big]_\m-\big[[x_1,x_2]_\h,y_3\big]=-\alpha([x_1,x_2]_\h)z_3.
	\end{equation}
	On the other hand, we have:
	\begin{equation}\label{Jac2}
		\big[[x_1,x_2]_\m,y_3\big]_\m=\mu_3\gamma y_3+\mu_3\eta z_3.
	\end{equation}
	Assume by contradiction that $\mu_3\neq 0$. Then, by identifying \eqref{Jac1} and \eqref{Jac2}, we deduce that:
	\begin{equation}\label{EE1}
		\gamma=0,\et \alpha([x_1,x_2]_\h)=-\mu_3\eta.
	\end{equation}
	A similar computation gives that:
	\begin{equation*}
		\big[[x_1,x_2]_\m,z_3\big]_\m=\big[[x_1,z_3]_\m,x_2\big]_\m+\big[x_1,[x_2,z_3]_\m\big]_\m-\big[[x_1,x_2]_\h,z_3\big]=\alpha([x_1,x_2]_\h)y_3.
	\end{equation*}
	Now, by a direct computation, we have:
	\begin{equation*}
		\big[[x_1,x_2]_\m,z_3\big]_\m=\mu_3\varepsilon y_3+\mu_3\xi z_3.
	\end{equation*}
	By identifying, we obtain that:
	\begin{equation}\label{EE2}
		\xi=0,\et \alpha([x_1,x_2]_\h)=\mu_3\varepsilon.
	\end{equation}
	In particular, combining \eqref{EE1} with \eqref{EE2}, we conclude that $\eta=-\varepsilon$. Moreover, using \eqref{Curv=0}, we get:
	\begin{equation*}
		0=\langle\R(y_3,z_3)x_1,x_2\rangle=\frac{\nu}{2}\big(\mu_2+\mu_3-\mu_1\big).
	\end{equation*}
	If $\mu_1=\mu_2+\mu_3$, then from \eqref{Lamd1andLamd2}, one can see that $\lambda_1=\lambda_2=0$, which is impossible. Thus, $\mu_1\neq \mu_2+\mu_3$, and therefore $\nu=0$. In summary, we have:
	\begin{equation*}
		[y_3,z_3]_\m=0,\quad [x_3,y_3]_\m=-\varepsilon z_3,\et [x_3,z_3]_\m=\varepsilon y_3.
	\end{equation*}
	By a straightforward computation, it follows that:
	\begin{equation}\label{Lamd3}
		\lambda_3=\ric(x_3,x_3)=\frac{(\mu_1+\mu_2)^2-\mu_3^2}{2}.
	\end{equation}
	One can see that the equalities \eqref{Lamd1andLamd2} and \eqref{Lamd3} are exactly the same as \eqref{la1dim3} and \eqref{la2la3dim3}. Hence, the remainder of the proof is identical to that of the case where $\rico:\m\to\m$ has four distinct eigenvalues. This completes the proof of our main theorem.
	
	\section{Appendix}\label{Sect5}
	
	In this appendix, we provide an elementary proof of the equivalence, for a Riemannian manifold, between the harmonicity of the curvature tensor and the Ricci-Codazzi condition \eqref{RicciCodazIntrod}. Let $(M,\met)$ be an $n$-dimensional Riemannian manifold and denote by $\n$ its Levi-Civita connection, i.e., the linear connection on $M$ given by the following Koszul's formula:
	\begin{eqnarray*}
		\langle\n_XY,Z \rangle&:=&\tfrac{1}{2}\Big\{X\left(\langle Y,Z\rangle\right)+Y\left(\langle Z,X\rangle\right)-Z\left(\langle X,Y\rangle\right)\\
		&&\qquad-\langle Y,[X,Z]\rangle-\langle Z,[Y,X]\rangle+\langle X,[Z,Y]\rangle\Big\},
	\end{eqnarray*}
	for any $X,Y,Z\in\vecf(M)$. It is well known that the Levi-Civita connection $\n$ is the unique torsion-free linear connection on $M$ which is compatible with the Riemannian metric $\met$ (see for instance, \cite[p. 122]{Lee}), meaning that:
	\begin{eqnarray}
		\mathcal{T}(X,Y)&:=&\n_XY-\n_YX-[X,Y]\,=\,0;\nonumber\\
		(\n\met)(X,Y,Z)&:=&X\left(\langle Y,Z\rangle\right)-\langle\n_XY,Z\rangle-\langle Y,\n_XZ\rangle\,=\,0.\nonumber
	\end{eqnarray}
	The curvature tensor of $\n$, denoted by $\mathcal{R}$, is the $(1,3)$-tensor field on $M$ defined by:
	\begin{equation*}
		\mathcal{R}(X,Y)Z:=\n_X\n_YZ-\n_Y\n_XZ-\n_{[X,Y]}Z.
	\end{equation*}
	Using the curvature tensor $\mathcal{R}$ and the Riemannian metric $\met$, one can easily obtain a $(0,4)$-tensor field on $M$, known as the Riemann curvature tensor of $\met$. More precisely, for any quadruple of vector fields $X,Y,Z,Z'\in\vecf(M)$, the Riemann curvature tensor of $\met$ is defined by:
	\begin{equation*}
		\underline{\mathcal{R}}(X,Y,Z,Z'):=\langle\mathcal{R}(X,Y)Z,Z' \rangle.
	\end{equation*}
	
	The next proposition summarizes all symmetries of the Riemann curvature tensor (cf. \cite[p. 79]{Petr}).
	
	\begin{proposition}
		For any $X,Y,Z,Z'\in\vecf(M)$, the following hold:
		\begin{enumerate}
			\item[$1.$] $\underline{\mathcal{R}}(X,Y,Z,Z')=-\underline{\mathcal{R}}(Y,X,Z,Z')$;\vspace{1mm}
			
			\item[$2.$] $\underline{\mathcal{R}}(X,Y,Z,Z')=-\underline{\mathcal{R}}(X,Y,Z',Z)$;\vspace{1mm}
			
			\item[$3.$] $\underline{\mathcal{R}}(X,Y,Z,Z')=\underline{\mathcal{R}}(Z,Z',X,Y)$;\vspace{1mm}
			
			\item[$4.$] Bianchi’s first identity: $\displaystyle\sum_{\circlearrowright (X,Y,Z)}\underline{\mathcal{R}}(X,Y,Z,Z')=0$, where
			\begin{equation*}
				\displaystyle\sum_{\circlearrowright (X,Y,Z)}\underline{\mathcal{R}}(X,Y,Z,Z'):=\underline{\mathcal{R}}(X,Y,Z,Z')+\underline{\mathcal{R}}(Z,X,Y,Z')+\underline{\mathcal{R}}(Y,Z,X,Z').
			\end{equation*}
		\end{enumerate}
	\end{proposition}
	
	Recall that, $\n$ induces linear connections in all tensor bundles over $M$ (e.g. \cite[p. 95]{Lee} and \cite[p. 57]{Petr}). In particular, for the curvature tensor $\mathcal{R}$, we have:
	\begin{eqnarray*}
		\left(\n_W\mathcal{R}\right)(X,Y,Z)&:=&\n_W\big(\mathcal{R}(X,Y)Z\big)-\mathcal{R}\left(\n_WX,Y\right)Z\\
		&&-\mathcal{R}\left(X,\n_WY\right)Z-\mathcal{R}\left(X,Y\right)\n_WZ.
	\end{eqnarray*}
	The curvature tensor $\mathcal{R}$ satisfies the Bianchi’s second identity (cf. \cite[p. 79]{Petr}):
	\begin{equation*}
		\displaystyle\sum_{\circlearrowright (W,X,Y)}\left(\n_W\mathcal{R}\right)(X,Y,Z)=0.
	\end{equation*}
	
	The following lemma establishes a relation between $\n\mathcal{R}$ and $\n\underline{\mathcal{R}}$.
	
	\begin{lemma}\label{LemmaCovDerRiemCurv}
		For any $X,Y,Z,Z',W\in\vecf(M)$, we have:
		\begin{equation*}
			\left(\n_W\underline{\mathcal{R}}\right)(X,Y,Z,Z')=\langle \left(\n_W\mathcal{R}\right)(X,Y,Z),Z' \rangle.
		\end{equation*}
	\end{lemma}
	
	\begin{proof}
		By definition,
		\begin{eqnarray*}
			\left(\n_W\underline{\mathcal{R}}\right)(X,Y,Z,Z')&=&W\Big(\underline{\mathcal{R}}(X,Y,Z,Z')\Big)-\underline{\mathcal{R}}\left(\n_WX,Y,Z,Z'\right)\\
			&&-\underline{\mathcal{R}}\left(X,\n_WY,Z,Z'\right)-\underline{\mathcal{R}}\left(X,Y,\n_WZ,Z'\right)\\
			&&-\underline{\mathcal{R}}\left(X,Y,Z,\n_WZ'\right).
		\end{eqnarray*}
		On the other hand, since $\n$ is compatible with $\met$, the following equality holds:
		\begin{equation*}
			W\Big(\underline{\mathcal{R}}(X,Y,Z,Z')\Big)=\langle \n_W\big(\mathcal{R}(X,Y)Z\big),Z'\rangle+\underline{\mathcal{R}}\left(X,Y,Z,\n_WZ'\right).
		\end{equation*}
		Thus, by substitution, we obtain:
		\begin{eqnarray*}
			\left(\n_W\underline{\mathcal{R}}\right)(X,Y,Z,Z')&=&\langle \n_W\big(\mathcal{R}(X,Y)Z\big),Z'\rangle-\langle\mathcal{R}\left(\n_WX,Y\right)Z,Z'\rangle\\
			&&-\langle\mathcal{R}\left(X,\n_WY\right)Z,Z'\rangle-\langle\mathcal{R}\left(X,Y\right)\n_WZ,Z'\rangle\\
			&=&\langle \left(\n_W\mathcal{R}\right)(X,Y,Z),Z' \rangle.
		\end{eqnarray*}
		Hence, completing the proof.
	\end{proof}
	
	Using all symmetries of the Riemann curvature tensor $\underline{\mathcal{R}}$, the Bianchi’s second identity, and Lemma \ref{LemmaCovDerRiemCurv}, one can easily deduce the following proposition.
	
	\begin{proposition}\label{symRiemCurv}
		For any $X,Y,Z,Z',W\in\vecf(M)$, we have:
		\begin{enumerate}
			\item[$1.$] $\left(\n_W\underline{\mathcal{R}}\right)(X,Y,Z,Z')=-\left(\n_W\underline{\mathcal{R}}\right)(Y,X,Z,Z')$;\vspace{1mm}
			
			\item[$2.$] $\left(\n_W\underline{\mathcal{R}}\right)(X,Y,Z,Z')=-\left(\n_W\underline{\mathcal{R}}\right)(X,Y,Z',Z)$;\vspace{1mm}
			
			\item[$3.$] $\left(\n_W\underline{\mathcal{R}}\right)(X,Y,Z,Z')=\left(\n_W\underline{\mathcal{R}}\right)(Z,Z',X,Y)$;\vspace{1mm}
			
			\item[$4.$] $\displaystyle\sum_{\circlearrowright (X,Y,Z)}\left(\n_W\underline{\mathcal{R}}\right)(X,Y,Z,Z')=0$;
			\item[$5.$] $\displaystyle\sum_{\circlearrowright (W,X,Y)}\left(\n_W\underline{\mathcal{R}}\right)(X,Y,Z,Z')=0$.
		\end{enumerate}
	\end{proposition}
	
	Recall the following well-known definition (see \cite{Abr,Bess,Der,Gry}):
	
	\begin{definition}
		An $n$-dimensional Riemannian manifold $(M,\met)$ is said to have harmonic curvature if the divergence of $\mathcal{R}$ vanishes, i.e.,
		\begin{equation*}
			\divr(X,Y):=\somi\left(\n_{E_i}\R\right)(E_i,X,Y)=0,
		\end{equation*}
		for any $X,Y\in\vecf(M)$, where $\{E_i\}_{1\leqslant i\leqslant n}$ is a local orthonormal frame for $(M,\met)$.
	\end{definition}
	
	\begin{remark}
		Note that, the divergence of the curvature $\divr$ is a well-defined $(1,2)$-tensor field on $M$. In fact, if $\{F_i\}_{1\leqslant i\leqslant n}$ is another local orthonormal frame for $(M,\met)$, then, on a small neighborhood, for any $i\in\{1,\ldots,n\}$ write $F_i=\sum_{j=1}^{n}f_{ji}E_j$. Since $\{E_i\}_{1\leqslant i\leqslant n}$ and $\{F_i\}_{1\leqslant i\leqslant n}$ are both orthonormal, it follows that the matrix $P:=(f_{ij})_{1\leqslant i,j\leqslant n}$ satisfies $P^{-1}=P^{\T}$. Thus, we compute:
		\begin{eqnarray*}
			\somi\left(\n_{F_i}\R\right)(F_i,X,Y)&=&\somijkn f_{ki}f_{ji}\left(\n_{E_k}\R\right)(E_j,X,Y)\\
			&=&\somijkn f_{ik}f_{jk}\left(\n_{E_i}\R\right)(E_j,X,Y)\\
			&=&\somijn (PP^{\T})_{ij}\left(\n_{E_i}\R\right)(E_j,X,Y)\\
			&=&\somi\left(\n_{E_i}\R\right)(E_i,X,Y).
		\end{eqnarray*}
		This proves that $\divr$ does not depend on the chosen local orthonormal frame. 
	\end{remark}
	
	The following lemma will be needed to establish the desired equivalence condition for the harmonicity of the curvature. Since, it was left as an exercise in many Riemannian geometry text books (see \cite[p. 83]{DoC}, \cite[p. 150]{Lee}, and \cite[p. 72]{Petr}), we will prove it here.
	
	\begin{lemma}\label{GeodFramLem}
		Let $(M,\met)$ be a Riemannian manifold of dimension $n$. For each $p\in M$, there is a neighborhood $U$ of $p$ and an othonormal frame $\{E_i\}_{1\leqslant i\leqslant n}$ defined on $U$ such that: 
		\begin{equation}\label{GeodFram}
			\n_{X_p}E_i=0,\qquad\forall\,X\in\vecf(U).
		\end{equation}
		Such a local orthonormal frame is called a geodesic frame.
	\end{lemma}
	
	\begin{proof}
		Take a normal coordinate chart $(U,x_i)$ centered at $p\in M$ (see \cite[p. 131]{Lee}). We know that there exists a positive real number $\varepsilon>0$ such that $\exp_p:\mathbf{B}_\varepsilon(0)\subset T_pM\to U$ is a diffeomorphism with $\exp_p(0)=p$. For any $q\in U$, consider the geodesic $\gamma:[0,1]\to U$ defined by $\gamma(t):=\exp_p(tv_q)$, where $v_q$ is the unique vector of $\mathbf{B}_\varepsilon(0)$ such that $\exp_p(v_q)=q$. Thus, for any $i\in\{1,\ldots,n\}$, we can define a vector field $E_i$ on $U$ as follows:
		\begin{equation*}
			E_{i|q}:=\tau_{01}^\gamma\left(\partial_{i|p}\right),
		\end{equation*}
		where $\{\partial_{i|p}\}_{1\leqslant i\leqslant n}$ is the orthonormal frame of $T_pM$ induced by the normal coordinate chart $(U,x_i)$, and $\tau_{01}^\gamma:T_pM\to T_qM$ is the parallel transport map along $\gamma$. It is clear, since $\tau_{01}^\gamma:T_pM\to T_qM$ is a linear isometry ($\n$ is compatible with $\met$), that $\{E_{i}\}_{1\leqslant i\leqslant n}$ is an orthonormal frame on $U$. So, to complete the proof, it only remains to check \eqref{GeodFram}. Let $\alpha_i:[0,1]\to U$ be the geodesic given by $\alpha_i(t):=\exp_p(t\tfrac{\varepsilon}{2}\partial_{i|p})$. Then, for any $j\in\{1,\ldots,n\}$, we have $E_{_j|\alpha_i(t)}=\tau_{0t}^{\alpha_i}\left(\partial_{j|p}\right)$. Thus, using the last formula in \cite[p. 109]{Lee}, we compute:
		\begin{eqnarray*}
			\n_{E_{i|p}}E_j&=&\tfrac{2}{\varepsilon}\n_{\tfrac{\varepsilon}{2}\partial_{i|p}}E_j\\
			&=&\frac{2}{\varepsilon}\dtz\left(\tau_{0t}^{\alpha_i}\right)^{-1}\Big(E_{_j|\alpha_i(t)}\Big)\\
			&=&\frac{2}{\varepsilon}\dtz\left(\tau_{0t}^{\alpha_i}\right)^{-1}\Big(\tau_{0t}^{\alpha_i}\left(\partial_{j|p}\right)\Big)\\
			&=&\frac{2}{\varepsilon}\dtz\partial_{j|p}\\
			&=&0.
		\end{eqnarray*}
		Hence, for any $X\in\vecf(U)$, write $X=\sum_{i=1}^{n}f_iE_i$, we have:
		\begin{equation*}
			\n_{X_p}E_i=\somi f_i(p)\n_{E_{i|p}}E_j=0.
		\end{equation*}
		This finishes the proof.
	\end{proof}

	Recall that, the Ricci tensor of $\n$, denoted by $\ric$, is the $(0,2)$-tensor field on $M$ given by:
	\begin{equation*}
		\ric(X,Y):=\tr\Big(Z\mapsto\R(X,Z)Y\Big).
	\end{equation*}
	The Ricci operator of $\n$, denoted by $\rico$, is the $(1,1)$-tensor field on $M$ defined by:
	\begin{equation*}
		\langle\rico(X),Y\rangle:=\ric(X,Y).
	\end{equation*}
	
	Analogously to the divergence of the curvature tensor, we define the divergence of the Riemann curvature tensor, denoted $\divm$, as follows:
	\begin{equation}\label{Div=0}
		\divm(X,Y,Z):=\somi\left(\n_{E_i}\underline{\mathcal{R}}\right)(E_i,X,Y,Z),
	\end{equation}
	for any $X,Y,Z\in\vecf(M)$, where $\{E_i\}_{1\leqslant i\leqslant n}$ is a local orthonormal frame for $(M,\met)$. By Lemma \ref{LemmaCovDerRiemCurv}, one can easily see that $\divm=0$ if and only if $\divr=0$.\vspace{2mm}
	
	The following theorem relates the Riemann curvature tensor with the Ricci tensor.
	
	\begin{theorem}\label{TheoCodaz}
		Let $(M,\met)$ be a Riemannian manifold of dimension $n$. The Ricci tensor $\ric$ and the divergence of the Riemann curvature tensor $\divm$ are related by:
		\begin{equation}\label{codazzi0}
			\divm(Z,X,Y)=\left(\n_X\ric\right)(Y,Z)-\left(\n_Y\ric\right)(X,Z),
		\end{equation}
		for all $X,Y,Z\in\vecf(M)$.
	\end{theorem}
	
	\begin{proof}
		First, to prove that \eqref{codazzi0} holds, it suffices to check it at an arbitrary point $p\in M$. Let $p\in M$, then by Lemma \ref{GeodFramLem}, let $\{E_i\}_{1\leqslant i\leqslant n}$ be a geodesic frame defined on a neighborhood of $p$. For $i\in\{1,\ldots,n\}$, using Lemma \ref{LemmaCovDerRiemCurv} and Proposition \ref{symRiemCurv}, we compute:
		\begin{eqnarray*}
			\left(\n_{E_i}\underline{\mathcal{R}}\right)(E_i,Z,X,Y)&=&-\left(\n_{E_i}\underline{\mathcal{R}}\right)(X,Y,Z,E_i)\\
			&=&\left(\n_{Y}\underline{\mathcal{R}}\right)(E_i,X,Z,E_i)+\left(\n_{X}\underline{\mathcal{R}}\right)(Y,E_i,Z,E_i)\\
			&=&\left(\n_{X}\underline{\mathcal{R}}\right)(Y,E_i,Z,E_i)-\left(\n_{Y}\underline{\mathcal{R}}\right)(X,E_i,Z,E_i)\\
			&=&\langle (\n_{X}\mathcal{R})\left(Y,E_i,Z\right),E_i\rangle-\langle (\n_{Y}\mathcal{R})\left(X,E_i,Z\right),E_i\rangle.
		\end{eqnarray*}
		On the other hand, we have:
		\begin{eqnarray*}
			\langle (\n_{X}\mathcal{R})\left(Y,E_i,Z\right),E_i\rangle&=&X\left(\langle\R(Y,E_i)Z,E_i\rangle\right)-\langle\R(Y,E_i)Z,\n_XE_i\rangle\\
			&&-\langle\R\left(\n_XY,E_i\right)Z,E_i\rangle-\langle\R\left(Y,\n_XE_i\right)Z,E_i\rangle\\
			&&-\langle\R\left(Y,E_i\right)\n_XZ,E_i\rangle\\
			&=&X\left(\langle\R(Y,E_i)Z,E_i\rangle\right)-\langle\R\left(\n_XY,E_i\right)Z,E_i\rangle\\
			&&-\langle\R\left(Y,E_i\right)\n_XZ,E_i\rangle-\langle\R(Y,E_i)Z,\n_XE_i\rangle\\
			&&-\langle\R\left(Y,\n_XE_i\right)Z,E_i\rangle.
		\end{eqnarray*}
		Thus, summing over $i$ and evaluating at $p$, we obtain that:
		\begin{eqnarray*}
			\somi\langle (\n_{X}\mathcal{R})\left(Y,E_i,Z\right),E_i\rangle_p&=&X_p\left(\somi\langle\R(Y,E_i)Z,E_i\rangle\right)-\somi\langle\R\left(\n_XY,E_i\right)Z,E_i\rangle_p\\
			&&-\somi\langle\R\left(Y,E_i\right)\n_XZ,E_i\rangle_p\\
			&=&X_p\left(\ric(Y,Z)\right)-\ric\left(\n_XY,Z\right)_p-\ric\left(Y,\n_XZ\right)_p\\
			&=&\left(\n_X\ric\right)(Y,Z)_p.
		\end{eqnarray*}
		By a similar argument, we obtain that:
		\begin{equation*}
			\somi\langle (\n_{Y}\mathcal{R})\left(X,E_i,Z\right),E_i\rangle_p=\left(\n_Y\ric\right)(X,Z)_p.
		\end{equation*}
		Hence, by substitution, we get:
		\begin{equation*}
			\divm(Z,X,Y)_p=\somi\left(\n_{E_i}\underline{\mathcal{R}}\right)(E_i,Z,X,Y)_p=\left(\n_X\ric\right)(Y,Z)_p-\left(\n_Y\ric\right)(X,Z)_p,
		\end{equation*}
		which proves \eqref{codazzi0}.
	\end{proof}
	
	As direct consequence of Theorem \ref{TheoCodaz}, we obtain the following:
	
	\begin{corollary}\label{CorolApx}
		A Riemannian manifold $(M,\met)$ has a harmonic curvature if and only if it is Ricci-Codazzi, meaning that its Ricci tensor satisfies \eqref{RicciCodazIntrod}.
	\end{corollary}

\vskip 0.2cm {\bf Acknowledgments.} The authors would like to express their sincere gratitude to Abdelhak Abouqateb (Cadi Ayyad University, Morocco) for his many insightful and constructive comments, which significantly improved the presentation of this paper. 

\end{document}